\documentclass[final,3p,times]{elsarticle}
\usepackage{amssymb}
\usepackage[colorlinks,
            linkcolor=blue,
            anchorcolor=blue,
            citecolor=blue]{hyperref}
\usepackage{graphicx} 
\usepackage{epstopdf}
\usepackage{subfigure}
\usepackage{booktabs} 
\usepackage{multirow}
\usepackage{CJK}
\usepackage{subfigure}
\usepackage{algorithm}
\usepackage{algorithmicx}
\usepackage{algpseudocode}
\usepackage{amsmath}
\usepackage{indentfirst}
\usepackage{caption}
\usepackage{times}
\usepackage{float}
\usepackage{stfloats}
\journal{Swarm and Evolutionary Computation}

\begin{document}
\begin{frontmatter}
\title{An enhanced POSTA based on Nelder-Mead simplex search and quadratic interpolation}
\tnotetext[1]{This study was driven by interest.}
\author{Tianyu Liu\corref{mycorrespondingauthor}}
\cortext[mycorrespondingauthor]{Corresponding author}
\ead{e1011072@u.nus.edu}

\begin{abstract}
State transition algorithm (STA) is a metaheuristic method for global optimization. Recently, a modified STA named parameter optimal state transition algorithm (POSTA) is proposed. In POSTA, the performance of expansion operator, rotation operator and axesion operator is optimized through a parameter selection mechanism. But due to the insufficient utilization of historical information, POSTA still suffers from slow convergence speed and low solution accuracy on specific problems. To make better use of the historical information, Nelder-Mead (NM) simplex search and quadratic interpolation (QI) are integrated into POSTA. The enhanced POSTA is tested against 14 benchmark functions with 20-D, 30-D and 50-D space. An experimental comparison with several competitive metaheuristic methods demonstrates the effectiveness of the proposed method.
\end{abstract}


\end{frontmatter}

\section{Introduction}
State transition algorithm (STA) \cite{STA2012,sta2020principle} is a novel metaheuristic method for global optimization and has been successfully applied in many fields \cite{zhou2019dynamic,wang2021hierarchical,zhou2021fast}. The basic idea of STA is to regard a solution as a state. Therefore the update of a solution can be considered as a state transition. In STA, four state transformation operators are designed to generate candidate solutions, namely expansion operator, rotation operator, axesion operator and translation operator. Among these operators, expansion operator is used for global search, rotation operator is for local search, axesion operator is for single-dimensional search and translation operator is designed for heuristic search \cite{POSTA2019}. Designed for different functionality, the state transformation operators are iteratively invoked to generate promising solutions until terminal conditions are satisfied.

Recently, a modified STA named parameter optimal state transition algorithm (POSTA) \cite{POSTA2019} is proposed. Based on a statistical study, POSTA provides a parameter selection mechanism to select optimal parameters for expansion operator, rotation operator and axesion operator. As a result, POSTA has better solution accuracy and stronger global exploration ability than STA. However, the historical information in POSTA is not sufficiently utilized. During the whole search process, the insufficient utilization of historical information leads to a slow convergence speed. In POSTA, the historical information is utilized by the translation operator. When a better solution emerges, the translation operator searches along the line between the current solution and the last solution with a maximal length of $\beta$. Obviously, the historical information is not sufficiently utilized, only two historical solutions are considered and only a linear transformation is applied for utilization. This inefficiency leads to a relatively slow convergence speed. In the later stage of the search, the lack of utilizing historical information leads to a low solution accuracy. As a metaheuristic method, POSTA is efficient to globally explore the whole solution space and localize the promising areas. But when the promising areas are reached, the lack of using historical information leads to a low probability of generating better solutions. This inefficiency leads to a relatively low solution accuracy. The above discussion implies that due to the insufficient utilization of historical information, POSTA still suffers from slow convergence speed and low solution accuracy on specific problems.

To speed up the convergence during the whole search process, Nelder-Mead (NM) simplex search \cite{nelder1965simplex,NM2018basic}
is introduced. NM simplex search is a classical direct search method with a fast convergence speed. For decades, combining NM simplex search with metaheuristic methods has been a popular way to speed up the convergence \cite{NM2021generalcombine,NMlocal2019swarm}. But none of the existing combination strategies have seen NM simplex search as a way to utilize the historical information. In this paper, a new combination strategy is proposed by applying NM simplex search to make use of the historical information.
In the proposed method, the NM simplex is used to store the historical information and the NM geometric transformations are applied to utilize the historical information. The NM simplex is a geometric figure consisting of $n+1$ vertices in n-dimensional space. Therefore, the NM simplex is able to store the information of $n+1$ historical solutions, which is more than the two historical solutions in translation operator. On the other hand, the NM geometric transformations consists of a sequence of distinct geometric transformations (reflection, contraction \dots), which are more comprehensive than the linear transformation in translation operator. Therefore, applying NM simplex search is feasible to utilize the historical information more sufficiently and speed up the convergence. To improve the solution accuracy in later stage of the search, a quadratic interpolation (QI) strategy is introduced. QI is a classical local search method and is widely used \cite{QI2021swarm,QI2020QISCA,QI2018,QI2019QIWOA}. In QI, three known points are used as the search agents to generate a quadratic curve. The generated curve is seen as an approximate shape of the target function. Therefore the extreme point of the quadratic function can be seen as an approximate optima of the target function. To make use of the historical information, three historical solutions are selected as the QI search agents to generate a new analytical solution with little computing cost. Therefore, QI is feasible to strengthen the exploitation capacity by utilizing historical information. By applying NM simplex search and QI in different stage of the search, the historical information of POSTA is utilized efficiently. Therefore, an enhanced POSTA based on NM simplex search and QI is proposed.

The main contributions of this paper can be summarized as follows: (1) A modified POSTA named the enhanced POSTA is proposed. (2) In the proposed method, the historical information is utilized more sufficiently to improve the performance. (3) The proposed method successfully combines the merits of the three distinct methods: the wide exploration capacity of POSTA, the fast convergence speed of NM simplex search and the strong exploitation capacity of QI. (4) When combining NM simplex search and POSTA, a new combination strategy that differs from the existing ones is proposed.

The remainder of the article is organized as follows.\ Section 2 reviews the algorithmic background of POSTA, NM simplex search and QI. In Section 3, the enhanced POSTA based on NM simplex search and QI is proposed.\ In Section 4, the efficiency and stability of the proposed method are illustrated through experimental comparisons.

\section{Background}
~\\
\textbf{2.1. Parameter optimal state transition algorithm (POSTA)}
~\

STA is a novel metahuristic method with applications in many fields \cite{zhou2020kernel, DONG2021, zhou2020using, zhou2020hybrid}. In STA, a solution can be regarded as a state, and the update of a solution can be seen as a state transition. In STA, based on state space representation, the general form of generating a solution is described as :

$$
\boldsymbol{s}_{k+1}=A_{k}\boldsymbol{s}_{k}+B_{k} \boldsymbol{u}_{k} \\
\eqno{(1)}
$$
where $\boldsymbol{s}_{k}$ and $\boldsymbol{s}_{k+1}$ are the current best state and the next state; $A_{k}$ and $B_{k}$ are state transition matrices that indicates the effect of state transformation operators; $\boldsymbol{u}_{k}$ is a function of $\boldsymbol{s}_{k}$ and historical states.

~\\
\textbf{A. State transformation operators}
~\

1)	Rotation transformation
$$
\boldsymbol{s}_{k+1}=\boldsymbol{s}_{k}+\alpha \frac{1}{n\left\|\boldsymbol{s}_{k}\right\|_{2}} R_{r} \boldsymbol{s}_{k}\eqno{(2)}
$$
where $\alpha$ is a positive parameter named rotation factor; $R_{r} \in \mathbb{R}^{n \times n}$ stands for a random matrix and all its elements are distributed in the range of [-1, 1]; $\|\cdot\|_{2}$ is the L2-norm of a vector. The rotation transformation is capable to generate candidate solutions in a hypersphere with a maximum radius of $\alpha$ and is designed for local search.

2) Translation transformation
$$
\boldsymbol{s}_{k+1}=\boldsymbol{s}_{k}+\beta R_{t} \frac{\boldsymbol{s}_{k}-\boldsymbol{s}_{k-1}}{\left\|\boldsymbol{s}_{k}-\boldsymbol{s}_{k-1}\right\|_{2}}
\eqno{(3)}$$
where $\beta$ is a positive constant named translation factor; $R_{t} \in \mathbb{R}$ stands for a random variable and all its elements belong to the range of [0,1]. Designed for heuristic search, the translation transformation is only invoked when a better solution is found by other operators.

3)	Expansion transformation
$$
\boldsymbol{s}_{k+1}=\boldsymbol{s}_{k}+\gamma R_{e}\boldsymbol{s}_{k}
\eqno{(4)}
$$
where $\gamma$ is a positive parameter named expansion factor;
$R_{e} \in \mathbb{R}^{n \times n}$  is a random diagonal matrix with its entries obeying the Gaussian distribution. The expansion transformation is designed to reinforce the capacity of global search and is able to search for solutions in the entire search space.

4)	Axesion transformation
$$
\boldsymbol{s}_{k+1}=\boldsymbol{s}_{k}+\delta R_{a}\boldsymbol{s}_{k}
\eqno{(5)}$$
where $\delta$ is a positive parameter named axesion factor; $R_{a} \in \mathbb{R}^{n \times n}$ is a random diagonal matrix with most of its elements obeying the Gaussian distribution and only a random one has a nonzero value. The axesion transformation is designed to strengthen the single-dimensional search.

~\\
\textbf{B. Parameter selection mechanism}
~\

According to a statistical study \cite{POSTA2019}, the parameter of the operators could be a crucial factor for the performance. To simplify the selection of the parameter, the candidate values for all parameters are taken from the set $\Omega=\{1,1 \mathrm{e}$-1, $1 \mathrm{e}$-2, $1\mathrm{e}$-3, $1\mathrm{e}$-4, $1 \mathrm{e}$-5, $1\mathrm{e}$-6, $1\mathrm{e}$-7, $1\mathrm{e}$-8\}. Then the value that leads to the optimal objective function value is chosen as the optimal value for the parameter. If the optimal parameter is denoted as $\tilde{a}^{*}$, the following formula represents the parameter selection:
$$
\tilde{a}^{*}=\underset{{\tilde{a}_{k} \in \Omega}}{\arg\min}f\left(\boldsymbol{s}_{k}+\tilde{a}_{k} \boldsymbol{\tilde{d}}_{k}\right)
\eqno{(6)}
$$
$$
\left.\left.\begin{array}{l}
\frac{1}{n\left\|\boldsymbol{s}_{k}\right\|_{2}} R_{r} \boldsymbol{s}_{k} \\
R_{e} \boldsymbol{s}_{k} \\
R_{a} \boldsymbol{s}_{k}
\end{array}\right\} \Rightarrow \boldsymbol{\tilde{d}}_{k}, \quad \begin{array}{c}
\alpha \\
\gamma \\
\delta
\end{array}\right\} \Rightarrow \tilde{a}_{k}.
\eqno{(7)}$$

~\\
\textbf{C. Algorithm procedure}
~\

To make use of the optimal parameter, the selected optimal value is kept for a period of time that is denoted as $T_{p}$. Specifically, the detailed procedures of the POSTA is as follows:

\begin{algorithmic}[1]
\State  $\textbf{repeat}$
\State {\ \ \ \ Best $\gets$ expansion\_w(objfun, Best, SE, $\Omega$)}
\State {\ \ \ \ Best $\gets$ rotation\_w(objfun, Best, SE, $\Omega$)}
\State {\ \ \ \ Best $\gets$ axesion\_w(objfun, Best, SE, $\Omega$)}
\State  {\textbf{until} certain criterion is  met}

\end{algorithmic}
where objfun indicates the objective function, SE indicates the number of candidate solutions in the candidate solution set and Best indicates the current best solution.

Meanwhile, rotation\_w($\cdot$) in the above pseudocode is given detailed explanations:

\begin{algorithmic}[1]
\State  { [Best,$\alpha$] $\gets$ update\_alpha(objfun, Best, SE, $\Omega$)}
\For {each $i\in [1,T_{p}]$}
    \State {Best $\gets$  rotation(objfun, Best, SE, $\Omega$)}
\EndFor
\end{algorithmic}
where rotation($\cdot$) represents the invocation of rotation operator and update\_alpha($\cdot$) represents the selection of optimal value for rotation factor. The parameter selection for expansion factor and axesion factor are similar to that for rotation factor. More details of POSTA can be found in \cite{POSTA2019}.

~\\
\textbf{2.2. Nelder-Mead simplex search}
~\\

Nelder-Mead simplex search is a popular direct search method \cite{ANMS2012} with a fast convergence speed. The idea of NM simplex search is to create a changeable simplex to approach and then approximate the optimum. For a D-dimensional function, it first initializes a simplex with $D+1$ vertices and iteratively replace the worst vertex by a better one. The new vertices in this process are generated by five geometric transformations, namely reflection, expansion, outside contraction, inside contraction and shrinkage. The detailed steps in NM simplex search are as follows:

\textbf{Step 0: Initialization.}
Starting from an original point $\mathbf{x}_{0}$ in D-dimensional space, the first simplex is generated by creating tiny perturbation in each dimension of $\mathbf{x}_{0}$ with a rate of 0.25\%. Specifically, the remaining $D$ vertices are generated as:

$$
 \mathbf{x}_{i}=\mathbf{x}_{0}+ \tau_{i}\mathbf{e}_{i}
 \eqno{(8)}$$
$$
 \tau_{i}=\left\{\begin{array}{ll}
0.05, & \text { if }\left(\mathbf{x}_{0}\right)_{i} \neq 0 \\
0.00025, & \text { if }\left(\mathbf{x}_{0}\right)_{i}=0
\end{array}\right.
\eqno{(9)}$$
where $i=1, \ldots, D+1$, $\left(\mathbf{x}_{0}\right)_{i}$ is the component of $\mathbf{x}_{0}$ in the i-th dimension and $\mathbf{e}_{i}$ is the unit vector in the i-th dimension.

\textbf{Step 1: Sorting.} The vertices $\mathbf{x}_{i}, \ i=1, \ldots, D+1$ are sorted and relabeled so that the function values are as:

$$
f\left(\mathbf{x}_{1}\right) \leq f\left(\mathbf{x}_{2}\right) \leq \ldots \leq f\left(\mathbf{x}_{D+1}\right)\eqno{(10)}
$$

\textbf{Step 2: Reflection.} Compute the reflection vertex $\mathbf{x}_{r}$  as:

$$
\mathbf{x}_{r}=\mathbf{x}_{c}+\eta\left(\mathbf{x}_{c}-\mathbf{x}_{D+1}\right)\eqno{(11)}
$$
where $\mathbf{x}_{c}=\left(\sum_{i=1}^{D} \mathbf{x}_{i}\right) /D $ and $\eta$ is the reflection coefficient. If  $f\left(\mathbf{x}_{1}\right) \leq f\left(\mathbf{x}_{r}\right)<f\left(\mathbf{x}_{D}\right)$, then replace $\mathbf{x}_{D+1}$\ with \ $\mathbf{x}_{r}$.

\textbf{Step 3: Expansion.}\ If $f\left(\mathbf{x}_{r}\right)<f\left(\mathbf{x}_{1}\right) $, compute the expansion vertex $\mathbf{x}_{e}$ as:

$$
\mathbf{x}_{e}=\mathbf{x}_{c}+\lambda\left(\mathbf{x}_{r}-\mathbf{x}_{c}\right)\eqno{(12)}
$$
where $\lambda$ is the expansion coefficient. If $f\left(\mathbf{x}_{e}\right)<f\left(\mathbf{x}_{r}\right)$, then replace $\mathbf{x}_{D+1}$ with $\mathbf{x}_{e}$, else with $\mathbf{x}_{r}$.

\textbf{Step 4: Outside Contraction(OC).}\ If $f\left(\mathbf{x}_{D}\right) \leq f\left(\mathbf{x}_{r}\right)<f\left(\mathbf{x}_{D+1}\right)$, then compute the OC vertex $\mathbf{x}_{\rm oc}$ as:

$$
\mathbf{x}_{oc}=\mathbf{x}_{c}+\mu\left(\mathbf{x}_{r}-\mathbf{x}_{c}\right)\eqno{(13)}
$$
where $\mu$ is the contraction coefficient. If $f\left(\mathbf{x}_{oc}\right) \leq f\left(\mathbf{x}_{r}\right)$, then replace\ $\mathbf{x}_{D+1}$ with $\mathbf{x}_{oc}$, else go to step 6.

\textbf{Step 5: Inside Contraction(IC).} If $ f\left(\mathbf{x}_{r}\right)\geq f\left(\mathbf{x}_{D+1}\right) $, then compute the IC vertex $\mathbf{x}_{ic}$ as:

$$
\mathbf{x}_{ic}=\mathbf{x}_{c}-\mu\left(\mathbf{x}_{r}-\mathbf{x}_{c}\right)\eqno{(14)}$$
where $\mu$ is the contraction coefficient. If $f\left(\mathbf{x}_{ic}\right)<f\left(\mathbf{x}_{D+1}\right)$ , then replace $\mathbf{x}_{D+1}$ with $\mathbf{x}_{ic}$, else go to step 6.

\textbf{Step 6: Shrinkage.} For $2 \leq i \leq D+1$:

$$
\mathbf{x}_{i}=\mathbf{x}_{1}+\nu\left(\mathbf{x}_{i}-\mathbf{x}_{1}\right)\eqno{(15)}
$$
where $\nu$ is the shrinkage coefficient.

\textbf{Step 7:}
If the termination condition is satisfied, stop, else go to step 1.
\begin{figure}[H]
\centering
\subfigure[Reflection] {
\label{nmgt1:a}
\includegraphics[width=0.2\columnwidth]{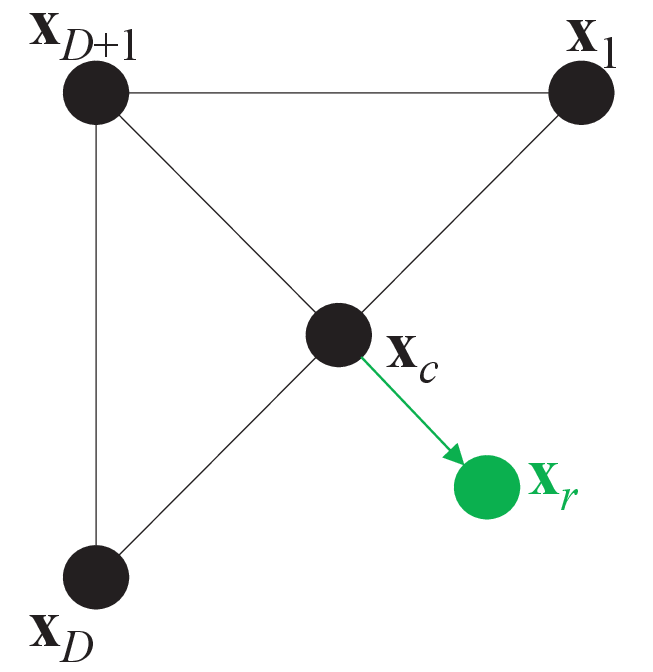}
}
\subfigure[Expansion] {
\label{nmgt1:b}
\includegraphics[width=0.22\columnwidth]{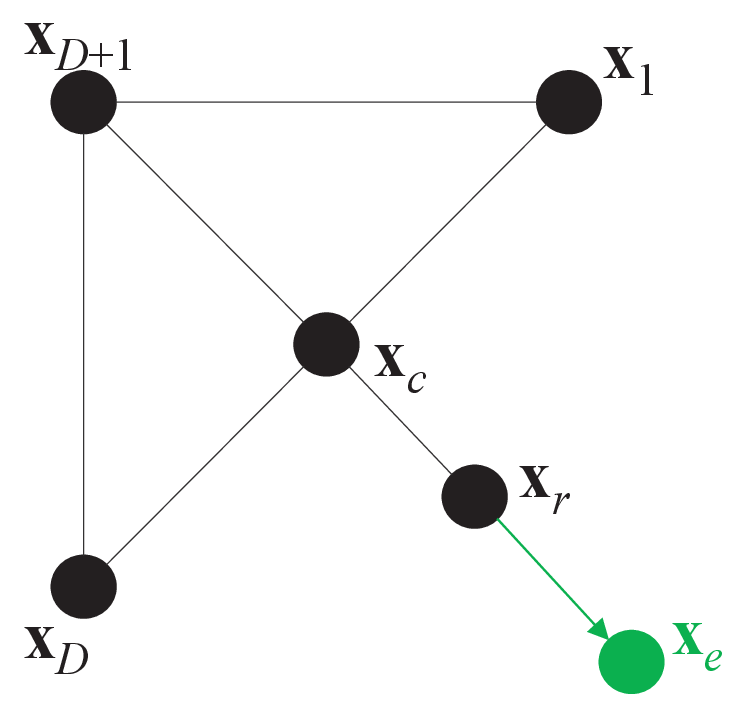}
}
\subfigure[Outside contraction] {
\label{nmgt1:c}
\includegraphics[width=0.21\columnwidth]{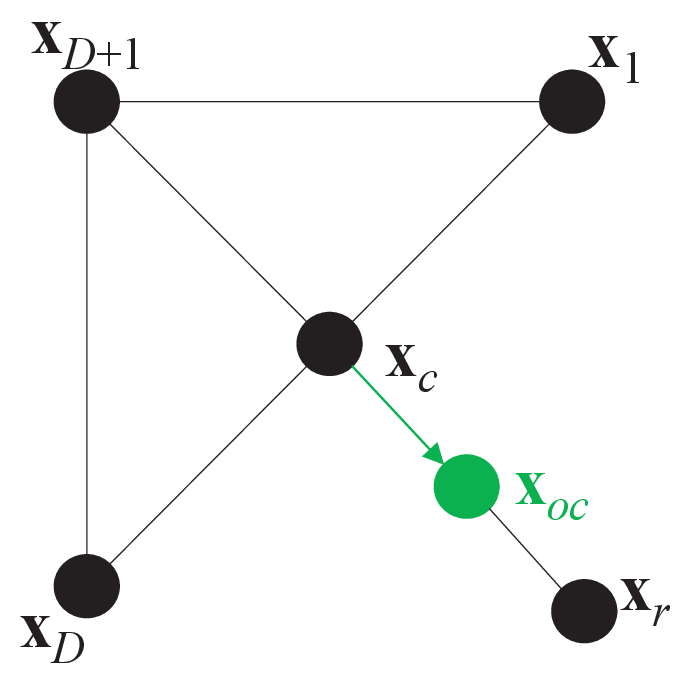}
}
\subfigure[Inside contraction] {
\label{nmgt1:d}
\includegraphics[width=0.19\columnwidth]{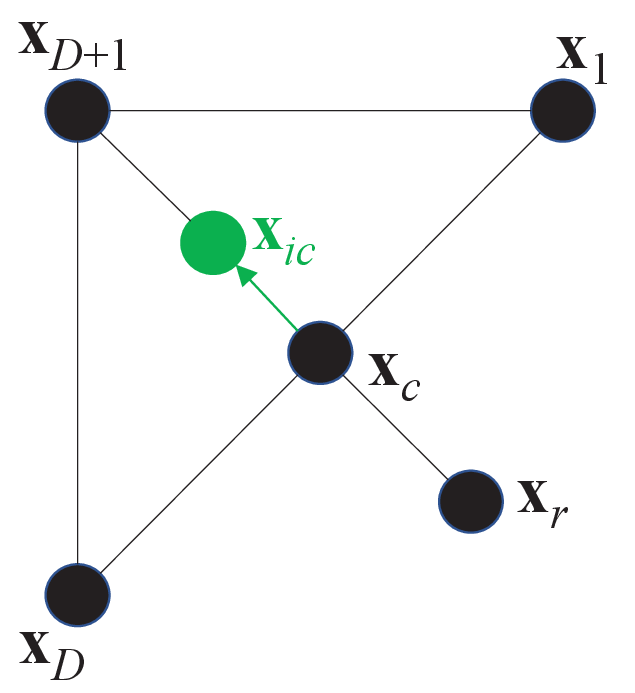}
}
\subfigure[Shrinkage] {
\label{nmgt1:e}
\includegraphics[width=0.2\columnwidth]{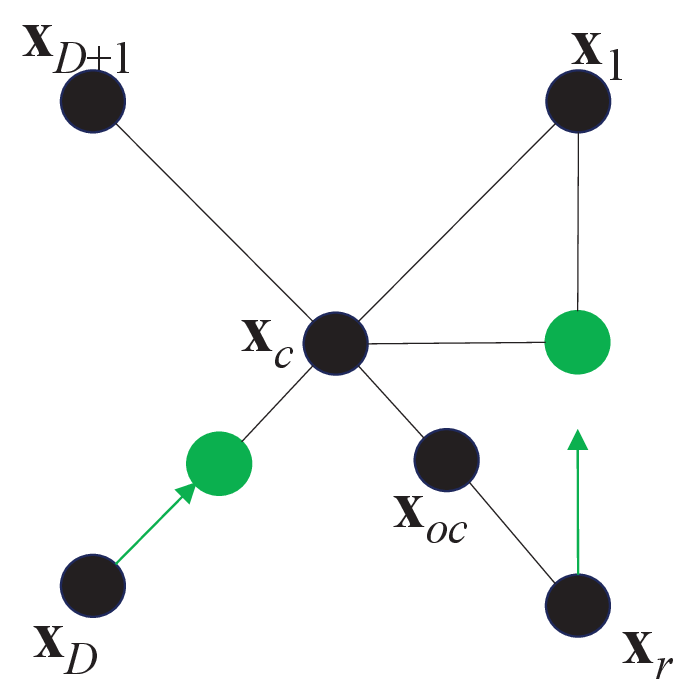}
}
\caption{The geometric transformations in NM simplex search.}
\label{nmgt}
\end{figure}

~\\
\textbf{2.3. Quadratic interpolation}
~\\

\begin{figure*}[h] 
\centering
\begin{equation}	
    x_{d}^{QI}=0.5 \frac{\left[\left(x_{d}^{best})^{2}-(x_{d}^{b})^{2}\right)\right] f\left(\mathbf{x}^{a}\right)+\left[\left(x_{d}^{a})^{2}-(x_{d}^{best})^{2}\right)\right] f\left(\mathbf{x}^{b}\right)+\left[\left(x_{d}^{b})^{2}-(x_{d}^{a})^{2}\right)\right] f\left(\mathbf{x}^{best}\right)}{\left(x_{d}^{best}-x_{d}^{b}\right) f\left(\mathbf{x}^{a}\right)+\left(x_{d}^{a}-x_{d}^{best}\right) f\left(\mathbf{x}^{b}\right)+\left(x_{d}^{b}-x_{d}^{a}\right) f\left(\mathbf{x}^{best}\right)},\ d=1,2, \ldots, D
    \tag {16}
    \label{QIequation2}
\end{equation}
\end{figure*}

QI is an analytical local search method which utilizes a parabola to fit the shape of the objective function. This method is capable to strengthen the exploitation capacity. In QI, three known points are used as the QI search agents to generate a quadratic curve, the generated quadratic curve is used as the approximate shape of the objective function, and the extreme point of the quadratic function can be used to approximate the optima of the objective function \cite{QI2017QIAEA}. Theoretically, when the the QI search agents are close enough to the global minimum, shape of the target function can be well approximated by the generated quadratic curve.

For a D-dimensional problem, three QI search agents $\mathbf{x}^{a}$, $\mathbf{x}^{b}$ and $\mathbf{x}^{best}$ are selected to generate a new solution $\mathbf{x}^{QI}$. Suppose that $$\mathbf{x}^{a}=\left(x_{1}^{a}, x_{2}^{a}, \ldots, x_{D}^{a}\right),$$
$$\mathbf{x}^{b}=\left(x_{1}^{b}, x_{2}^{b}, \ldots, x_{D}^{b}\right),$$ $$\mathbf{x}^{best}=\left(x_{1}^{best}, x_{2}^{best}, \ldots, x_{D}^{best}\right)$$ are the three distinct search agents, then the new solution $$\mathbf{x}^{QI}=\left(x_{1}^{QI}, x_{2}^{QI}, \ldots, x_{D}^{QI}\right)$$ can be calculated as Eq. (\ref{QIequation2}), where $f\left(\mathbf{x}^{a}\right)$, $f\left(\mathbf{x}^{b}\right)$ and $f\left(\mathbf{x}^{best}\right)$ are the fitness values of the three search agents respectively. As a strong local search method, QI is capable to strengthen the exploitation capacity but could easily lead to a local minimum. Therefore, QI should only be invoked in later stage of the search.

\section{Proposed method}
~\\
\textbf{3.1. NM-POSTA}
~\\

\begin{figure*}[h]

    \centering
        \subfigure[Translation operator] {
        \includegraphics[width=0.4\columnwidth]{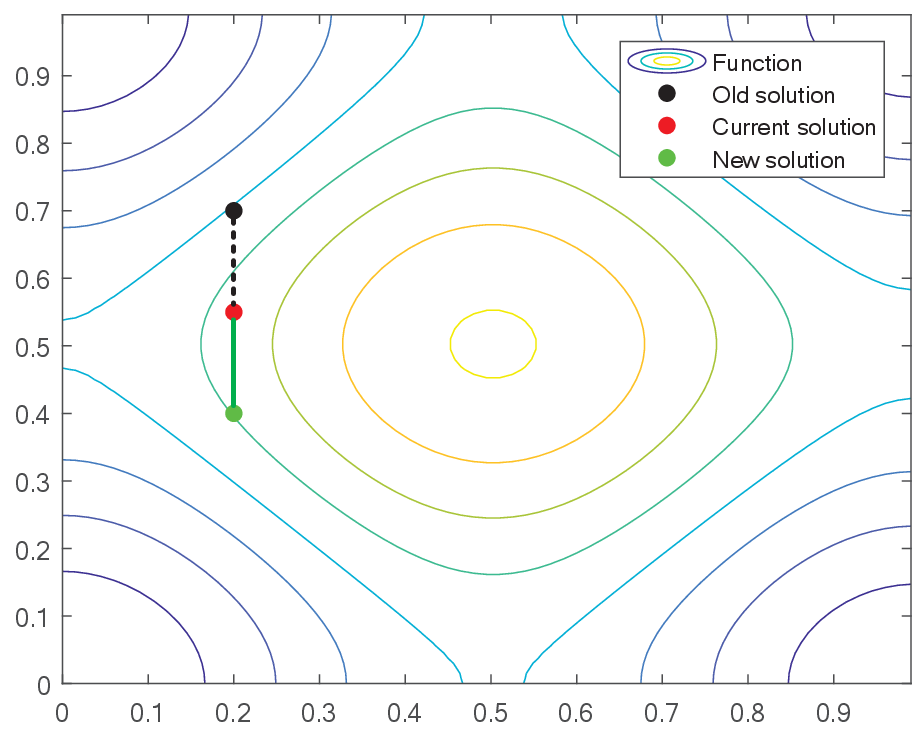}
        \label{utilize-comparison:a}
    }
        \subfigure[The proposed mechanism] {
        \includegraphics[width=0.4\columnwidth]{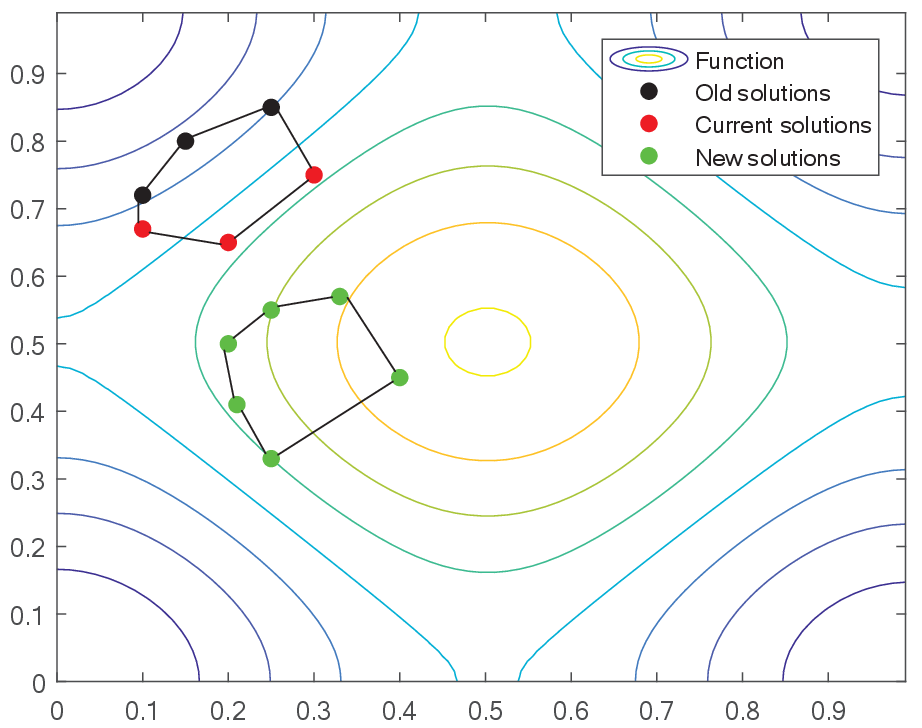}
        \label{utilize-comparison:b}
    }
    \caption{Comparison in utilizing history information.}
    \label{utilize-comparison}
\end{figure*}

\begin{figure*}[h]
\centering
\subfigure[Translation operator] {
 \label{trnm:a}
\includegraphics[width=0.6\columnwidth]{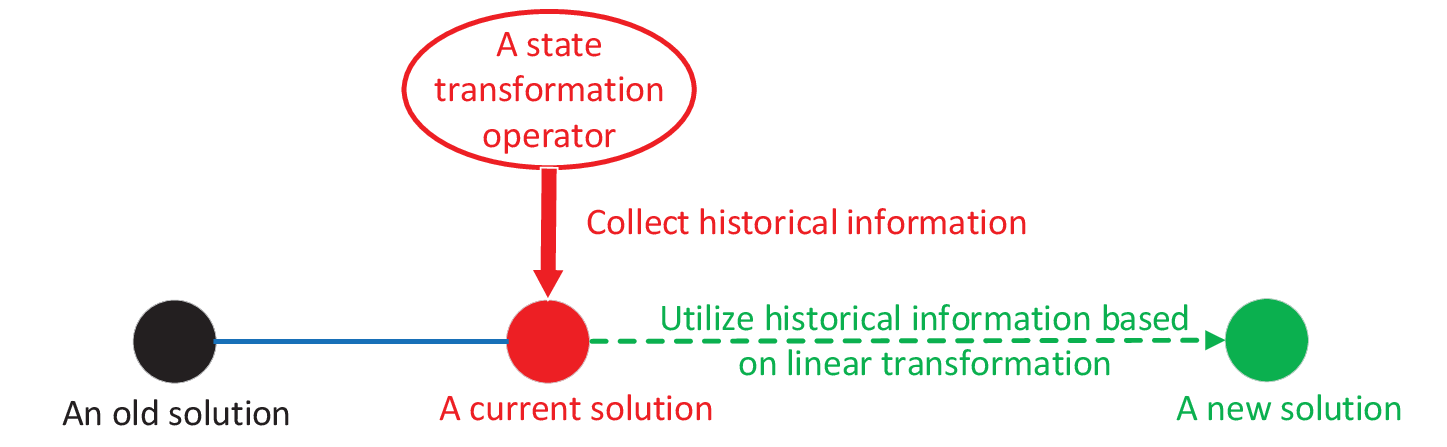}
}
\subfigure[The proposed mechanism] {
\label{trnm:b}
\includegraphics[width=0.9\columnwidth]{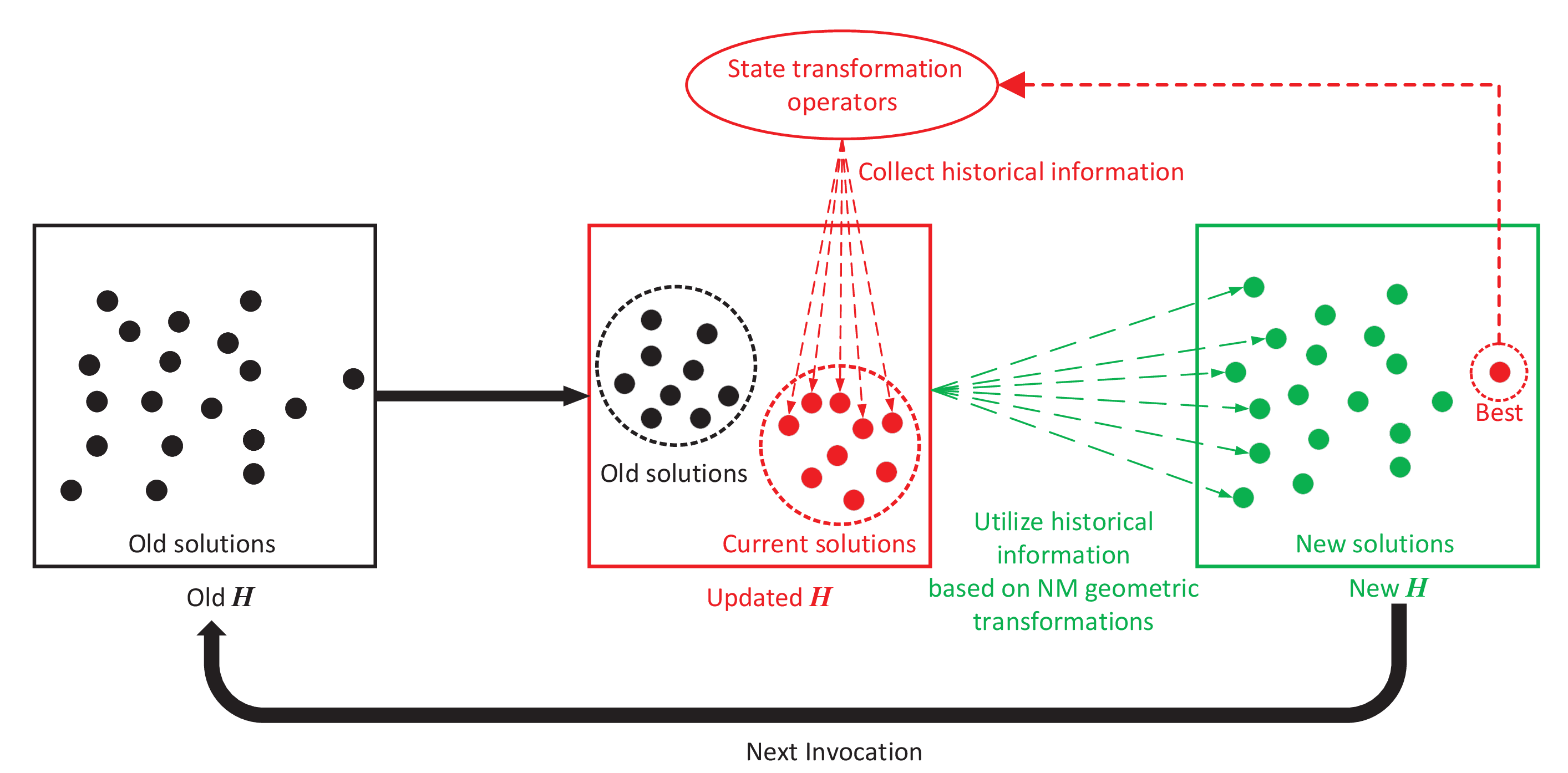}
}
\caption{Comparison in collecting and utilizing historical information. }
\label{trnm}
\end{figure*}

In POSTA, the historical information is not utilized sufficiently. As shown in Fig. \ref{utilize-comparison:a} and Fig. \ref{trnm:a} , the historical information is utilized by translation operator. However, only the latest two historical solutions are collected and only a linear search is applied to utilize the information. This inefficiency leads to the low convergence speed of POSTA on specific problems.

In the past decades, combining NM simplex search with metaheuristic methods has been a popular way to speed up the convergence. The general combination strategies between NM simplex search and metaheuristic methods can be summarized into two types: the staged pipelining type combination and the eagle strategy type combination \cite{NMABC2009}. In the first type, NM simplex search is applied to the superior individuals in the population \cite{NMWOA2019}. In the second type, NM simplex search is applied in the exploitation stage as a local search method \cite{xu2019EAGLE,NMlocal2019swarm}. But none of the existing combinations have seen NM simplex search as a way to utilize the information of historical solutions.

In this section, a POSTA enhanced with NM simplex search (NM-POSTA) is proposed. In NM-POSTA, a historical information mechanism based on NM simplex search is proposed to make use of the historical information. In the proposed mechanism, the historical information is collected based on a collection strategy, stored in the NM simplex, and then utilized by the NM geometric transformations. The detailed historical information mechanism in NM-POSTA is illustrated as follows.

~\\
\textbf{A. The historical information set}
~\

In NM-POSTA, a historical information set called $\textbf{\emph{H}}$ is used to store the historical solutions generated by the state transformation operators. $\textbf{\emph{H}}$ has the capacity to store the coordinates of $D+1$ vertices in D-dimensional space. Therefore $\textbf{\emph{H}}$ can be seen as an NM simplex. For a D-dimensional problem, $\textbf{\emph{H}}$ is capable to store the information of $D+1$ historical solutions.

The historical information contained in $\textbf{\emph{H}}$ is more sufficient than that in the translation operator. As shown in Fig. \ref{utilize-comparison} and Fig. \ref{trnm}, $\textbf{\emph{H}}$ contains information about a set of old solutions and a set of current solutions, but only an old solution and a current solution are considered in translation operator.

~\\
\textbf{B. Utilization of the historical information}
~\

\floatname{algorithm}{Algorithm}
\renewcommand{\algorithmicrequire}{\textbf{Input:}}
\renewcommand{\algorithmicensure}{\textbf{Output:}}
\begin{algorithm}
\caption{Utilize historical information based on NM geometric transformations.}  \label{hist-utilization}
\begin{algorithmic}[1] 
\Require Updated historical information set $\textbf{\emph{H}} $
\Ensure  New historical information set $\textbf{\emph{H}}$, $Best$
\State {initial simplex $\gets$ $\textbf{\emph{H}}$}
\For{each $i\in [1,D+1]$}
\State{NM geometric transformations}
\EndFor
\State {$\textbf{\emph{H}}$ $\gets$ new \ simplex}
\State {$Best$ $\gets$ best solution in $\textbf{\emph{H}}$}
\end{algorithmic}
\end{algorithm}

In NM-POSTA, the NM geometric transformations are applied for utilization. As shown in Fig. \ref{nmgt}, the NM geometric transformations consists of a sequence of distinct geometric transformations, which are more comprehensive than the linear transformation in translation operator. As a result, NM-POSTA is able to utilize the historical information more comprehensively, as shown in Fig. \ref{utilize-comparison} and Fig. \ref{trnm}.

The detailed steps of utilizing historical information in NM-POSTA is illustrated in Algorithm \ref{hist-utilization}. When Algorithm \ref{hist-utilization} is invoked, $\textbf{\emph{H}}$ is input as the initial simplex. After that, the NM geometric transformations are applied to the initialized simplex. In the inner iterations, the NM geometric transformations are run for $D+1$ times.

Considering that each iteration of NM geometric transformations produces a new vertex approximately, running $D+1$ times approximately produces a new simplex with $D+1$ new vertices. Therefore a new simplex with $D+1$ new vertices is generated, and is then stored in $\textbf{\emph{H}}$. As shown in Fig. \ref{trnm:b}, the new $\textbf{\emph{H}}$ contains $D+1$ new solutions and the top solution in $\textbf{\emph{H}}$ is denoted as $Best$. $Best$ is then sent to the state transformation operators as the input of current best solution.

~\\
\textbf{C. Collection of historical information}
~\

An appropriate collection strategy of history information is proposed to provide more promising input for utilization. Overall, the historical information is collected based on the collection strategy, stored in $\textbf{\emph{H}}$ and utilized by the NM geometric transformations. The utilization of historical information in NM-POSTA is illustrated in Algorithm \ref{hist-utilization}. If an invocation of utilization is terminated, then the historical information needs to be re-collected before next utilization. Therefore, between the invocations of utilization, a collection strategy for historical information is considered.

In the collection strategy, two type of solutions are considered: old solutions and current solutions. As shown in \ref{trnm:b}, old solutions are the solutions generated in the last invocation of utilization. After the termination of the last utilization, the state transformation operators start to generate solutions from $Best$. Between the invocations of utilization, those solutions generated by POSTA are considered as current solutions.

This collection strategy is seen as an expanded version of that in the translation operator. As shown in Fig. \ref{utilize-comparison:a} and Fig. \ref{trnm:a}, translation operator uses an old solution and a current solution as the historical information. In NM-POSTA, a set of old solutions and a set of current solutions are used as the historical information, as shown in Fig. \ref{utilize-comparison:b} and Fig. \ref{trnm:b}.

\floatname{algorithm}{Algorithm}
\renewcommand{\algorithmicrequire}{\textbf{Input:}}
\renewcommand{\algorithmicensure}{\textbf{Output:}}
\begin{algorithm}
\caption{Collect historical information}
\label{hist-collection}
\begin{algorithmic}[1] 
\Require Old historical information set $\textbf{\emph{H}}$, current  solution $\mathbf{x}_{current}$
\Ensure  Updated historical information set $\textbf{\emph{H}}$
\State {replace the worst solution in $\textbf{\emph{H}}$ with $\mathbf{x}_{current}$}
\State {\textbf{return} updated $\textbf{\emph{H}}$}
\end{algorithmic}
\end{algorithm}

\begin{figure}[t]
\centering
\includegraphics[width=0.55\columnwidth]{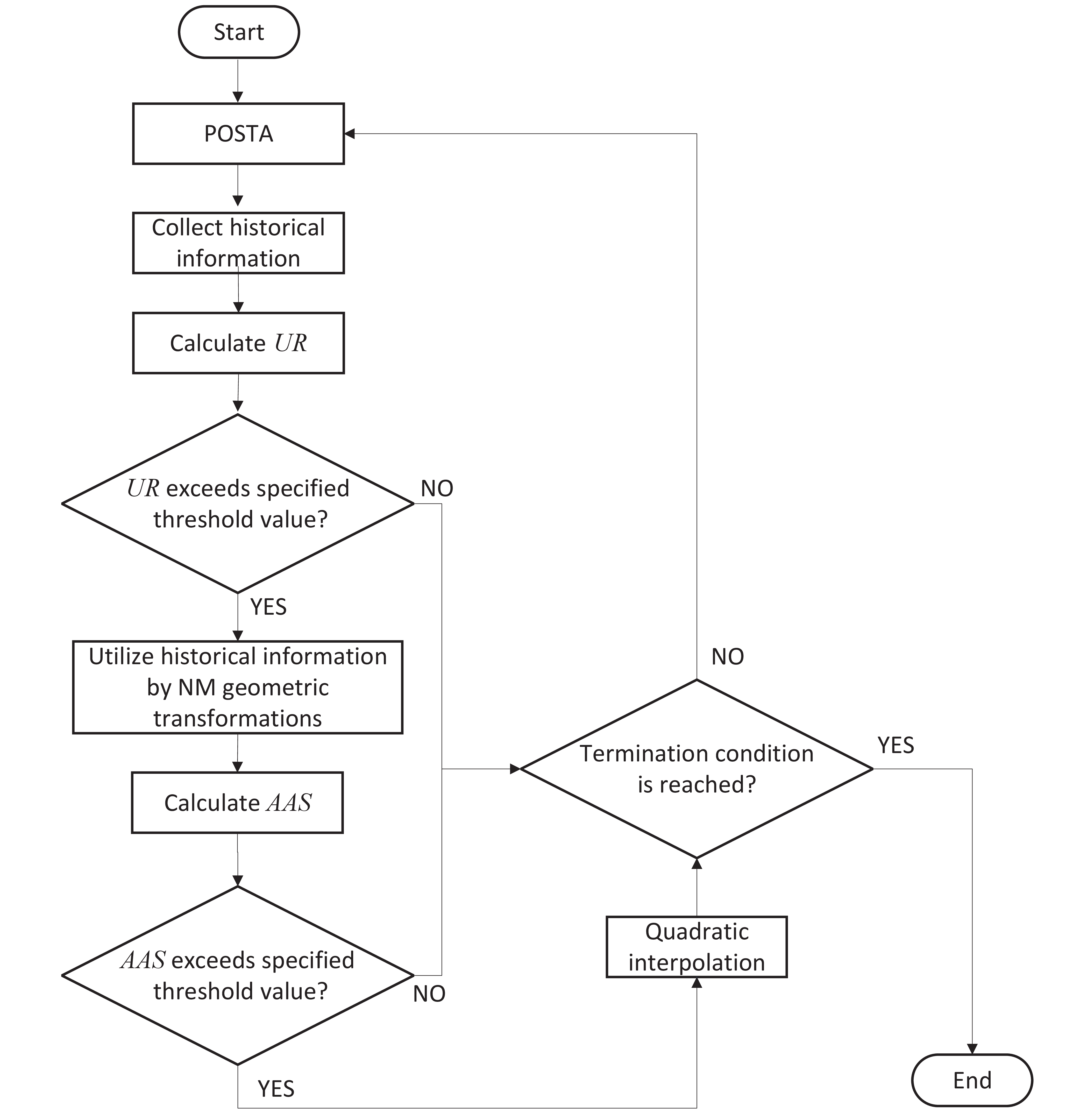}
\caption{Flow chart of the enhanced POSTA.}
\label{nmqiflowchart}
\end{figure}

The collection strategy is implemented by updating $\textbf{\emph{H}}$. As shown in Algorithm \ref{hist-collection}, a current solution $\mathbf{x}_{current}$ is used to update $\textbf{\emph{H}}$ in an invocation. In the implementation, Algorithm \ref{hist-collection} is invoked iteratively to update $\textbf{\emph{H}}$ with multiple current solutions. As shown in Fig. \ref{trnm:b}, $\textbf{\emph{H}}$ is updated by multiple current solutions before been utilized. To control the ratio between old solutions and current solutions in  $\textbf{\emph{H}}$, a parameter named update rate (UR) is proposed:

$$UR=\frac{n_{cs}}{n_{os}+n_{cs}}\eqno{(17)}$$
where $n_{os}$ and $n_{cs}$ are namely the number of old solutions and the number of current solutions in $\textbf{\emph{H}}$. As shown in Fig. \ref{nmqiflowchart}, $UR$ is calculated after each update, and only when $UR$ exceeds certain threshold value will the utilization be invoked.

~\\
\textbf{3.2. Properties of NM-POSTA}
~\\

\begin{figure*}[h]
\centering
\subfigure[Solution path of POSTA] {
\label{path:a}
\includegraphics[width=0.4\columnwidth]{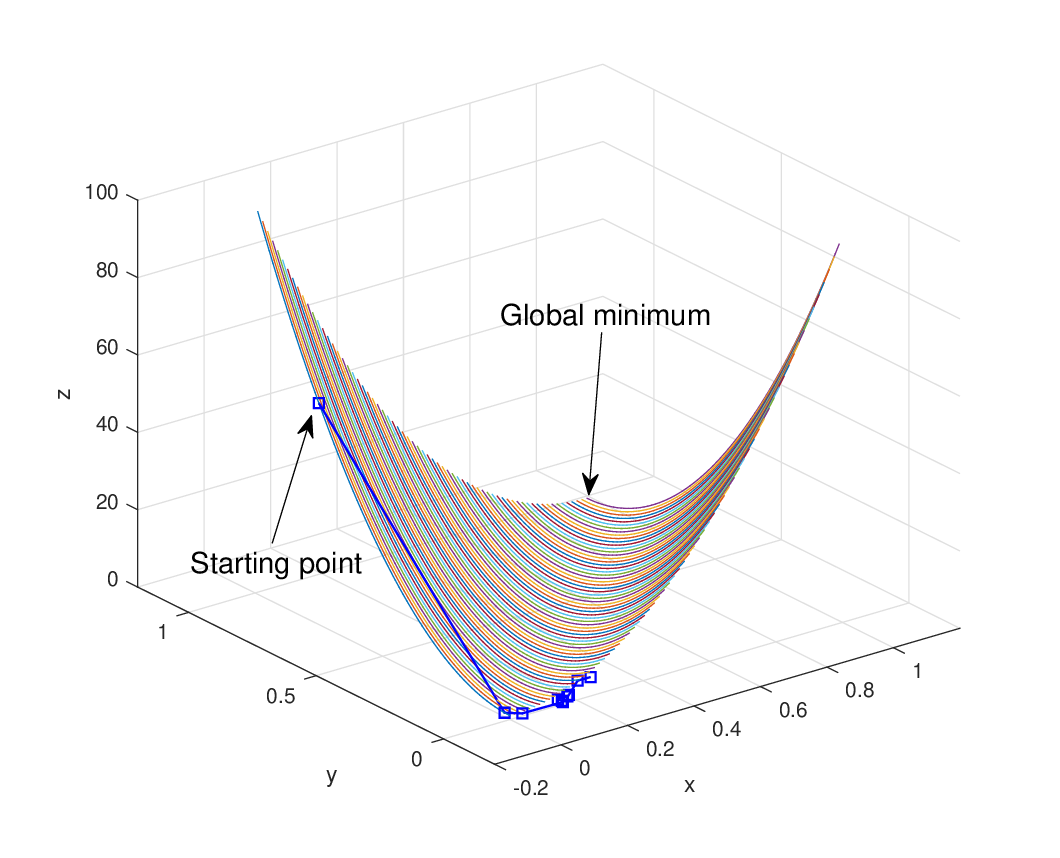}
}
\subfigure[Solution path of NM-POSTA] {
 \label{path:b}
\includegraphics[width=0.4\columnwidth]{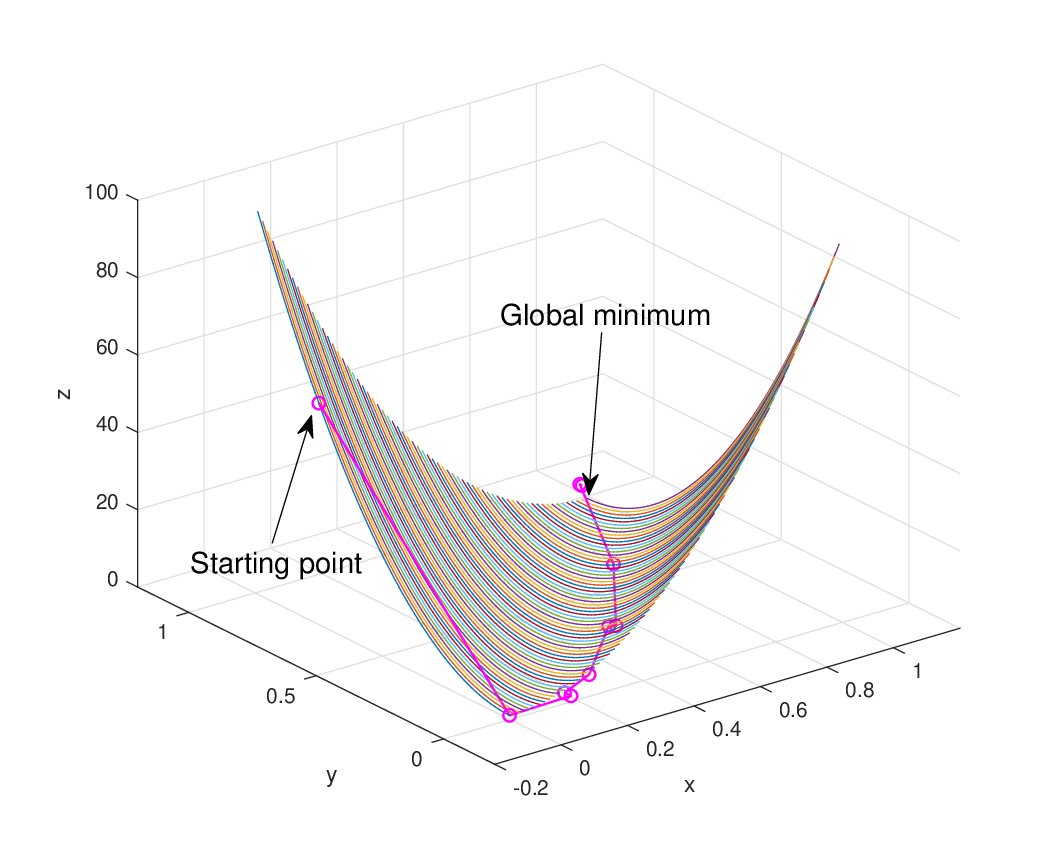}
}
\caption{Solution paths of POSTA and NM-POSTA on Rosenbrock function. }
\label{path}
\end{figure*}

In this section, some properties of NM-POSTA are briefly illustrated by a comparison between NM-POSTA and POSTA. To illustrate the effect of utilizing historical information, two classical test functions are used, namely Rosenbrock (F3) and Rastrigin (F7) from Table \ref{table3}. NM-POSTA and POSTA are tested against the two functions in 2-D space and each method is only terminated when certain criterion is satisfied. After termination, the corresponding number of function evaluations (FEs) is recorded. For both POSTA and NM-POSTA, $SE$ is set to 50 and $T_{p}$ is set to 10. In NM-POSTA, threshold value for $UR$ is set to 0.5 and parameters in NM simplex search are set as: $\{\eta,\ \lambda,\ \mu,\ \nu\}=\{1,\ 2,\ 0.5,\ 0.5\}$, which are the same as that in the standard implementation of NM simplex search \cite{nelder1965simplex}.

~\\
\textbf{A. Faster convergence speed}
~\

\begin{table*}[h]
\centering
\setlength{\belowcaptionskip}{10pt}
\caption{ \textbf{Table 1: }Statistical results for POSTA and NM-POSTA with D = 2. }
\renewcommand\arraystretch{1.15}
\scalebox{1}{
    \begin{tabular}{p{1.5cm}p{2.5cm}p{2.5cm}p{0.1cm}p{2.5cm}p{2.5cm}cccccc}
    \toprule
    \multirow{2}[0]{*}{function} & \multicolumn{2}{c}{POSTA}& & \multicolumn{2}{c}{NM-POSTA}  \\
    \cmidrule{2-3}\cmidrule{5-6}\
          & Ave FEs  & Success rate& & Ave FEs  & Success rate       \\
    \midrule
    F3    & 1.08E+04 & 30/30& & 4.54E+03 & 30/30   \\
    F7    & 1.51E+03 & 30/30& & 1.47E+03 & 30/30   \\
\bottomrule
\end{tabular}}%
\label{table1}%
\end{table*}%

To demonstrate the faster convergence speed of NM-POSTA, a 2-D Rosenbrock test function is used. This unimodal function has a global minimum (1,1) that lies in a narrow, parabolic valley. The valley is easy to reach but further convergence is very difficult \cite{Rosenbrock}.\ In the test for Rosenbrock function, the methods are terminated if the current solution $\mathbf{x}_{current}$ satisfies:
$$\left|{f}(\mathbf{x}_{current})-{f}({Best}^{*})\right| \leq \epsilon\eqno{(18)}$$
where ${Best}^{*}$ indicates the global minimum. If the specified accuracy $\epsilon$ is met, the test is considered a `success'. In this section, $\epsilon$ is set at 1e-8.

The statistical results given in Table \ref{table1} reveal that NM-POSTA is able to reach the same accuracy with much less FEs. In Fig. \ref{path}, the solution paths of POSTA and NM-POSTA on the Rosenbrock function are portrayed, where both methods start from the starting point (0, 0.75) and only the first 10 solutions are plotted in the figures for the convenience of observation. It is obvious that the solution path of NM-POSTA is much more efficient. By comparing the two paths, it can be seen that the both methods quickly reach the valley. After that, POSTA searches aimlessly in the valley and moves slowly while NM-POSTA finds the promising direction and moves towards the optimum efficiently. This brief experiment demonstrates that the utilization of historical information can lead to a more promising search direction, therefore NM-POSTA has a faster convergence speed than POSTA.

~\\
\textbf{B. Consistent global stability }
~\

To demonstrate the stability on avoiding local minimums, a 2-D Rastrigin test function is used. This function is a multimodal function with many local minimums.

In the test for Rastrigin function, the methods are only terminated when the global optimum is found. If the global minimum is found, the test is considered a `success'. The statistical results given in Table \ref{table1} show that NM-POSTA gets the same success rate with slightly less FEs. It is obvious from the results that the NM-POSTA remains stable in avoiding local minimums.

~\\
\textbf{3.3. The enhanced POSTA (NMQI-POSTA)}
~\\

POSTA, NM simplex search and QI are methods with distinct characteristics. As a novel metaheuristic method, POSTA has strong global exploration capacity but lacks convergence speed and exploitation capacity. As a classical direct search method, NM simplex search has a fast convergence speed but is not stable in global search. QI is able to approximate the global minimum by an analytical solution but is only effective when the search agents are very close to the global minimum.

To combine the merits of the three methods, NM simplex search and QI are applied to utilize the historical information of POSTA. Therefore, an enhanced POSTA based on NM simplex search and QI (NMQI-POSTA) is proposed. As shown in Fig. \ref{nmqiflowchart}, the proposed method uses an eagle strategy \cite{eagleorigin} to maximize the efficiency and stability. In the exploration stage, NM simplex search utilizes the historical information of POSTA to generate promising solutions. As a result, the convergence speed is continuously accelerated. In the exploitation stage, QI utilizes the historical information to improve the local search capacity. As a result, the solution accuracy is increased.

To invoke QI in later stage of the search, the average accuracy of the solutions in the historical information set $\textbf{\emph{H}}$ is calculated in each iteration as:

$$\small
AAS=\left|\frac{f\left(\textbf{x}_{1}\right) + f\left(\textbf{x}_{2}\right) + \ldots + f\left(\textbf{x}_{D+1}\right)}{D+1} -{f}({Best}^{*})\right|\eqno{(19)}
$$
where ${Best}^{*}$ indicates the global best solution.

If $AAS$ meets the certain threshold value, the solutions are considered to be in the neighborhood of the global optimum. Then QI will be invoked to increase the solution accuracy. In the implementation of QI, two QI search agents $\textbf{x}^{a}$ and $\textbf{x}^{b}$ are randomly chosen from the historical information set $\textbf{\emph{H}}$ and the third search agent $\textbf{x}^{best}$ is the current best solution.

For comparison, a POSTA enhanced with QI (QI-POSTA) is also implemented. The strategy of QI-POSTA is basically the same as NMQI-POSTA. The only difference is that the NM simplex search is disabled in QI-POSTA.

\section{Experimental results and discussions}

In this section, the proposed method is tested against a set of benchmark functions. To guarantee fairness and avoid arbitrariness, each test in the following experiments is performed 30 times independently. All methods are implemented in MATLAB R2017b environment with an AMD R7-4800H, 2.90 GHz processor on a 64-bit Windows 10 operating system.

\begin{table*}[h]
\centering
\setlength{\belowcaptionskip}{12pt}
\caption{ {\textbf{Table 2: }}Benchmark function set.}
\renewcommand\arraystretch{2}
\scalebox{0.7}{
\begin{tabular}{p{3cm}p{8cm}p{2cm}p{1cm}p{1.5cm}
lllll}
\toprule
   Funtion name    & Equation & Range & ${f}_{\min }$ & Modality  \\
  \midrule

  Elliptic &
  $\mathrm{F}1=\sum_{i=2}^{D}\left(10^{6}\right)^{(i-1) /(D-1)} \cdot x_{i}^{2}$ & $[-100,100]$  & 0 & Unimodal   \\

  Penalized 1 &

  $\begin{array}{c}
    \mathrm{F}2=\frac{\pi}{D}\left\{\sum_{t=1}^{D-1}\left(y_{i}-1\right)^{2}\left[1+\sin \left(\pi y_{i+1}\right)\right]+\left(y_{D}-1\right)^{2}+\right. \\
    \left.\left(10 \sin ^{2}\left(\pi y_{1}\right)\right)\right\}+\sum_{i=1}^{D} u\left(x_{i}, 10,100,4\right) \\
    y_{i}=1+\frac{x_{1}+1}{4} \\
    u\left(x_{i}, a, k, m\right)=\left\{\begin{array}{ll}
    k\left(x_{i}-a\right)^{m}, & x_{i}>a \\
    0, & -a \leq x_{i} \leq a \\
    k\left(-x_{i}-a\right)^{m}, & x_{i}<-a
    \end{array}\right.
    \end{array} $ & $[-50,50]$  & 0  & Multimodal \\

  Rosenbrock &
  $\mathrm{F}3=\sum_{i=1}^{D-1}\left[100\left(x_{i+1}-x_{i}^{2}\right)^{2}+\left(x_{i}-1\right)^{2}\right] $& $[-30,30]$  & 0 & Unimodal \\

  Schwefel 1.2 & $\mathrm{F}4=\sum_{i=1}^{D}\left(\sum_{j=1}^{1} x_{j}\right)^{2} $ & [-100,100]  & 0 & Unimodal \\

  Schwefel 2.4 & $\mathrm{F}5=\sum_{i=1}^{D}\left[\left(x_{i}-1\right)^{2}+\left(x_{1}-x_{i}^{2}\right)^{2}\right] $ & [0,10]  & 0 & Multimodal \\

  Sphere &
  $\mathrm{F}6=\sum_{i=1}^{D} x_{i}^{2} $ & [-100,100]  & 0 & Unimodal \\

  Rastrigin &
  $\mathrm{F}7=\sum_{i=1}^{D}\left[x_{i}^{2}-10 \cos \left(2 \pi x_{i}\right)+10\right] $ & [-5.12,5.12]  & 0 & Multimodal \\

  Griewank &
  $\mathrm{F}8=\sum_{i=1}^{D} \frac{x_{i}^{2}}{4000}-\prod_{i=1}^{D} \cos \left(x_{i} / \sqrt{i}\right)+1$ & [-60,60]  & 0 & Multimodal \\

  Sum squares &
  $\mathrm{F}9=\sum_{i=2}^{D} i x_{i}^{2}$ & [-10,10]  & 0 & Unimodal \\

  Levy and Montalvo 1 &
  $\begin{array}{c}
    \mathrm{F}10=\frac{\pi}{D}\left(10 \sin ^{2}\left(\pi y_{1}\right)+\sum_{i=1}^{D-1}\left(y_{i}-1\right)^{2}[1+\right. \\
    \left.\left.10 \sin ^{2}\left(\pi y_{i+1}\right)\right]+\left(y_{D}-1\right)^{2}\right), \quad y_{i}=1+\frac{1}{4}\left(x_{i}+1\right)
  \end{array}$ & [-10,10]  & 0 & Multimodal \\

  Zakharov &
  $\mathrm{F}11=\sum_{i=1}^{D} x_{i}^{2}+\left(\sum_{i=1}^{D} 0.5 i x_{i}\right)^{2}+\left(\sum_{i=1}^{D} 0.5 i x_{i}\right)^{4}$ & [-5,10]  & 0 & Unimodal \\

  Schwefel 2.2.2 &
  $\mathrm{F}12=\sum_{i=1}^{D}\left|x_{i}\right|+\prod_{i=1}^{D}\left|x_{i}\right|$ & [-10,10]  & 0 & Unimodal \\

  Cigar &
  $\mathrm{F}13=x_{1}^{2}+10^{6} \sum_{i-2}^{D} x_{i}^{6}$ & [-100,100]  & 0 & Unimodal \\

  Csendes &
  $\mathrm{F}14=\sum_{i=1}^{D} x_{i}^{6}\left(2+\sin \frac{1}{x_{i}}\right)$ & [-1,1]  & 0 & Multimodal \\
\bottomrule
\end{tabular}}%
\label{table3}%
\end{table*}%

In Section 4.2, the POSTA families (POSTA, NM-POSTA, QI-POSTA and NMQI-POSTA) are comprehensively tested against the benchmark functions. In Section 4.3, the enhanced POSTA is compared with several metaheuristic methods on the benchmark functions.

~\\
\textbf{4.1. Benchmark functions and parameter settings}
~\\

The details of the benchmark functions are shown in Table \ref{table3}. Notably, ${f}_{\min}$ and $Range$ indicate the global optimal value of the objective function and the boundary of the search space. Generally, in Table \ref{table3}, there are two types of benchmark functions, called unimodal and multimodal. A function with a single minimum in the specified range is called unimodal. Unlike unimodal, multimodal functions have many local minimums that the method may be trapped in.

For all members in the POSTA family, $SE$ is set to 50 and $T_{p}$ is set to 10. In NM-POSTA, threshold value for $UR$ is set to 0.5. In QI-POSTA and NMQI-POSTA, threshold value for $AAS$ is set to 1e-6. To guarantee fairness for the comparison in Section 4.3, the population sizes in the other metaheuristic methods are all set to 50, which are the same as the value of parameter $SE$ in the enhanced POSTA.

~\\
\textbf{4.2. Comparison among the POSTA families}
~\\
\begin{figure*}[h]
\centering
\subfigure[F1] {
\label{NM2:a}
\includegraphics[width=0.2\columnwidth]{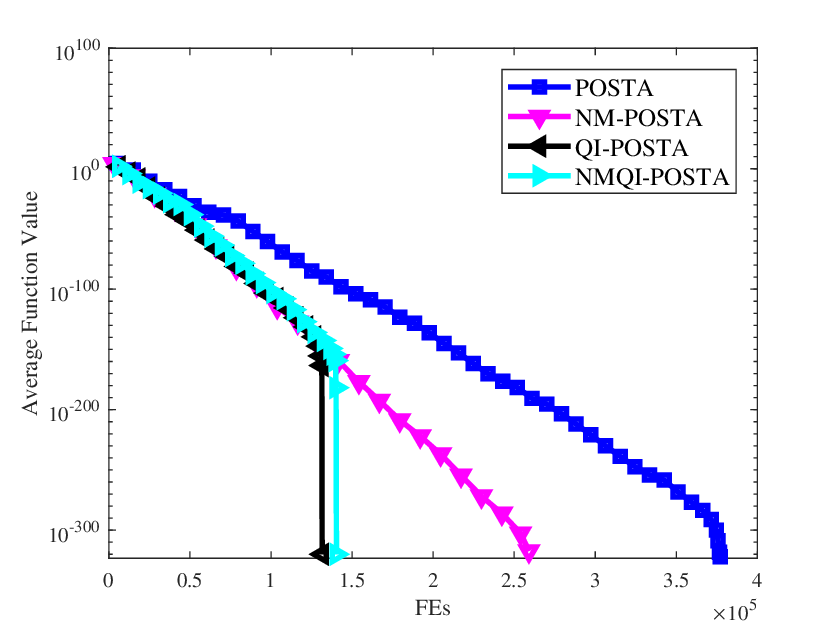}
}
\subfigure[F2] {
\label{NM2:b}
\includegraphics[width=0.2\columnwidth]{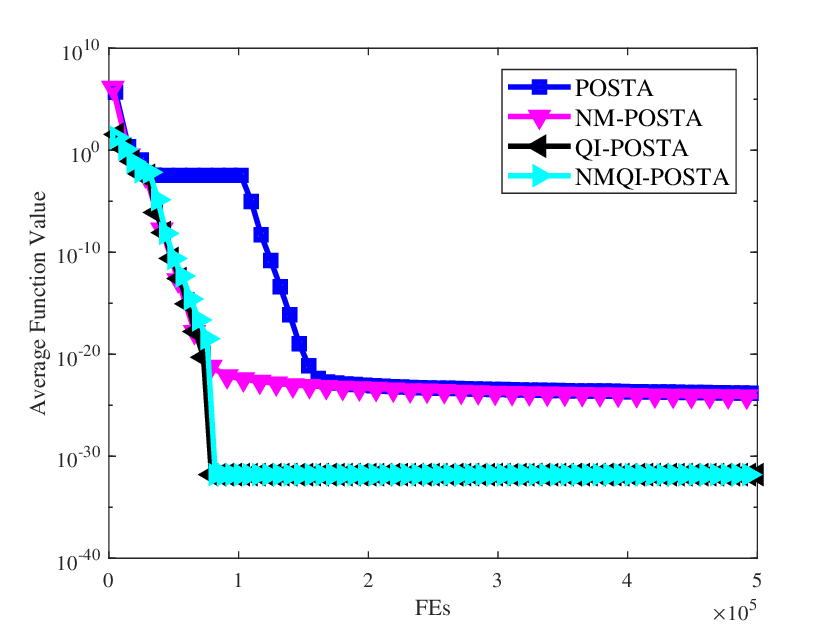}
}
\subfigure[F3] {
\label{NM2:c}
\includegraphics[width=0.2\columnwidth]{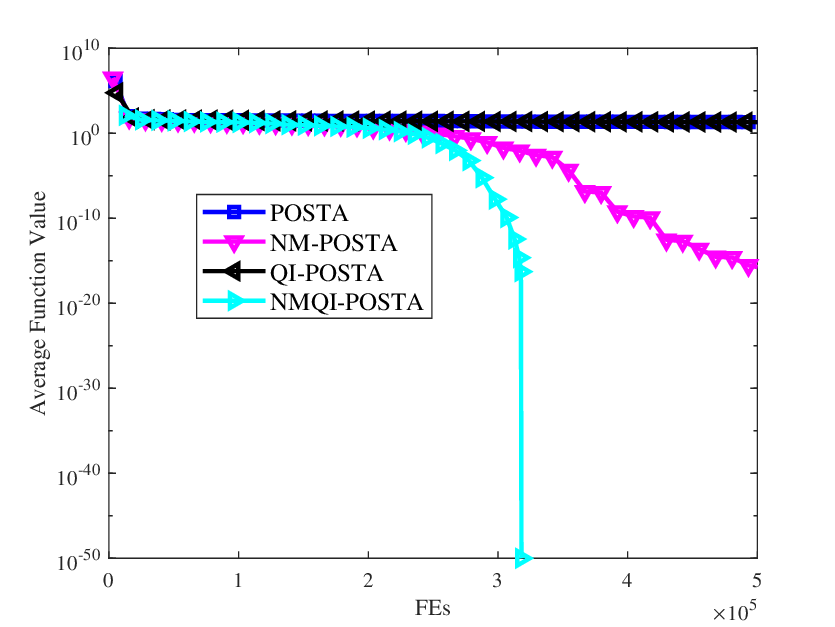}
}
\subfigure[F4] {
 \label{NM2:d}
\includegraphics[width=0.2\columnwidth]{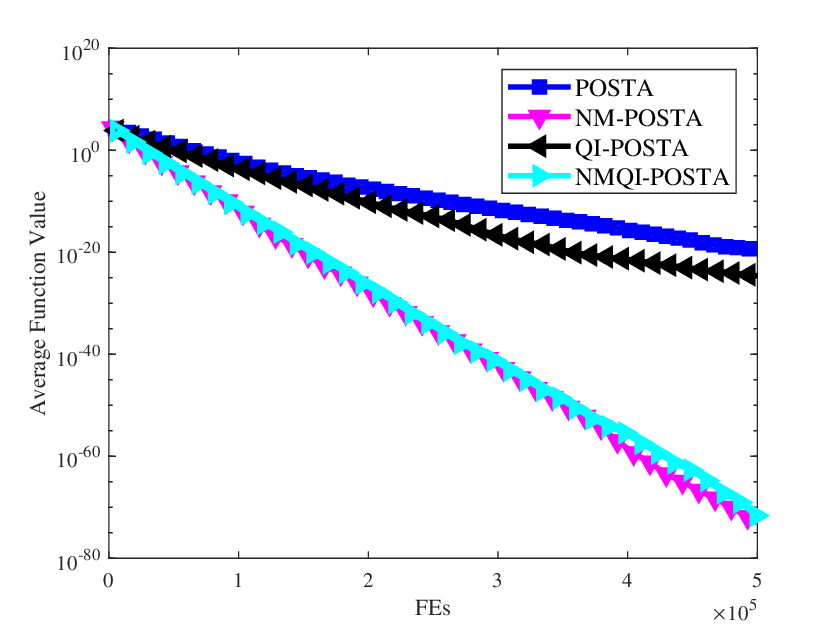}
}
\subfigure[F5] {
\label{NM2:e}
\includegraphics[width=0.2\columnwidth]{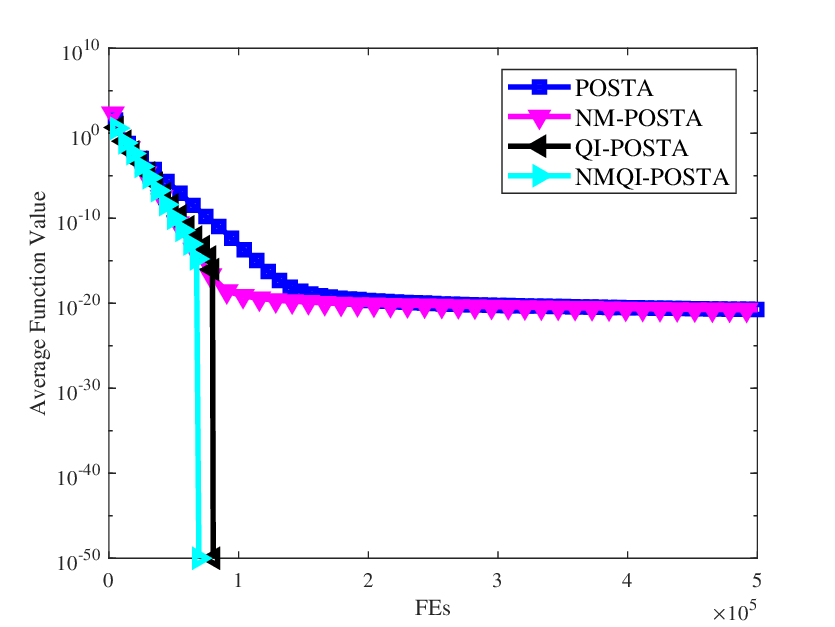}
}
\subfigure[F6] {
\label{NM2:f}
\includegraphics[width=0.2\columnwidth]{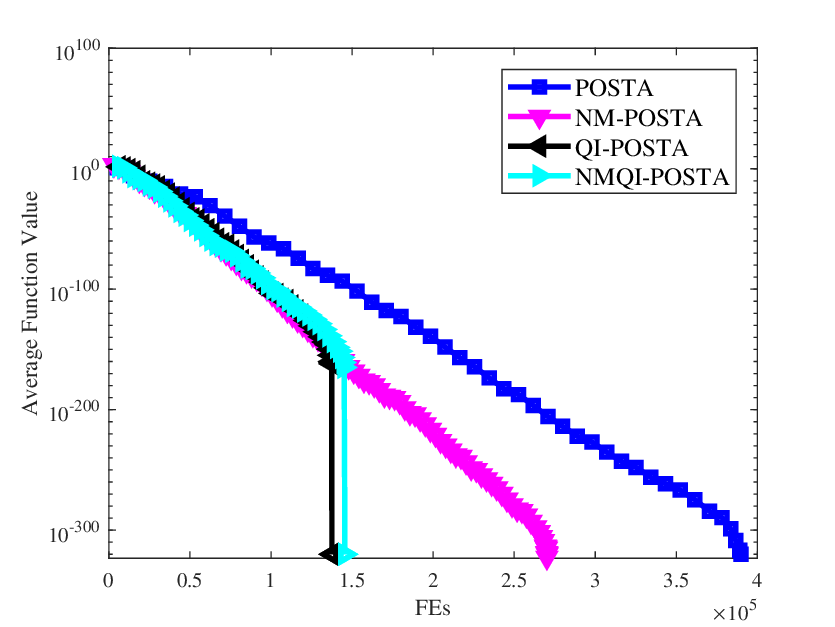}
}
\subfigure[F7] {
\label{NM2:g}
\includegraphics[width=0.2\columnwidth]{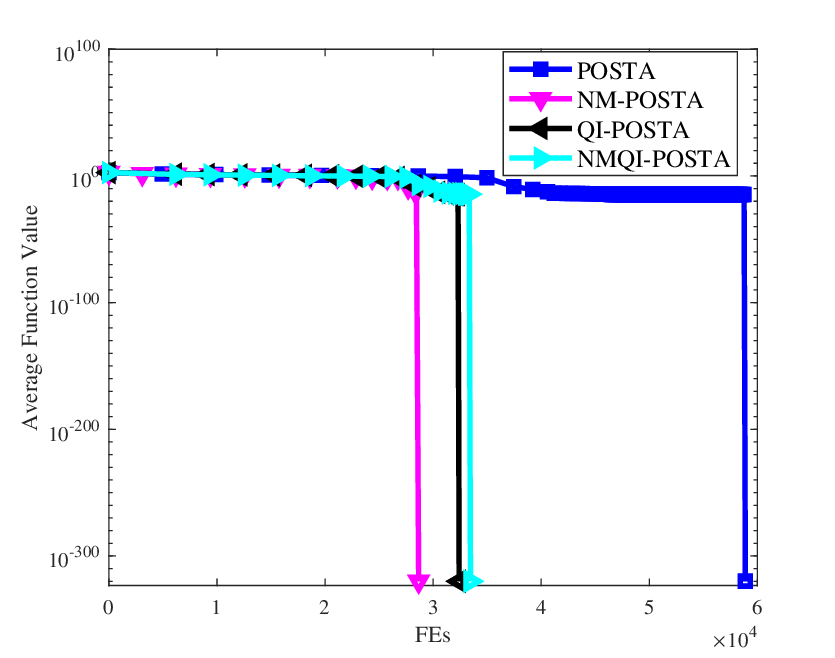}
}
\subfigure[F8] {
\label{NM2:h}
\includegraphics[width=0.2\columnwidth]{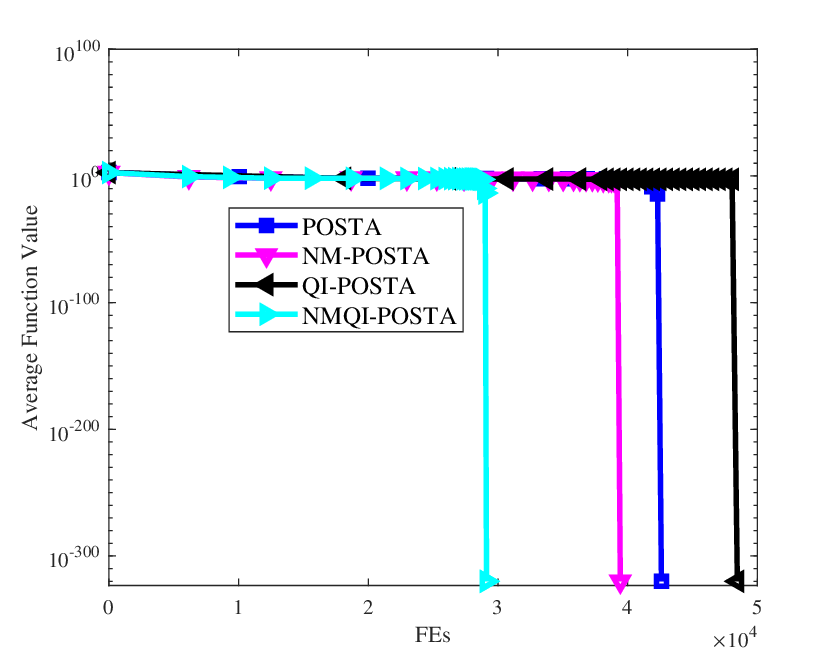}
}
\subfigure[F9] {
\label{NM2:i}
\includegraphics[width=0.2\columnwidth]{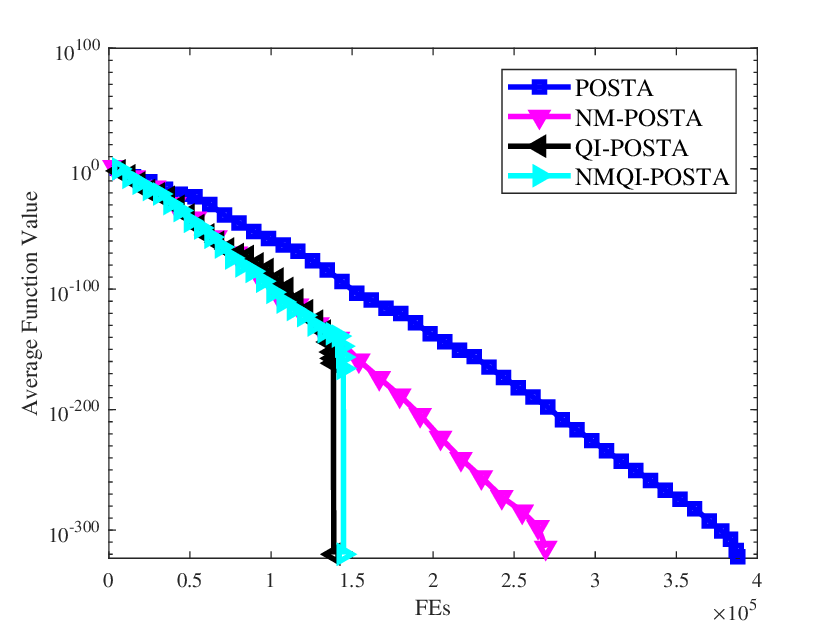}
}
\subfigure[F10] {
\label{NM2:j}
\includegraphics[width=0.2\columnwidth]{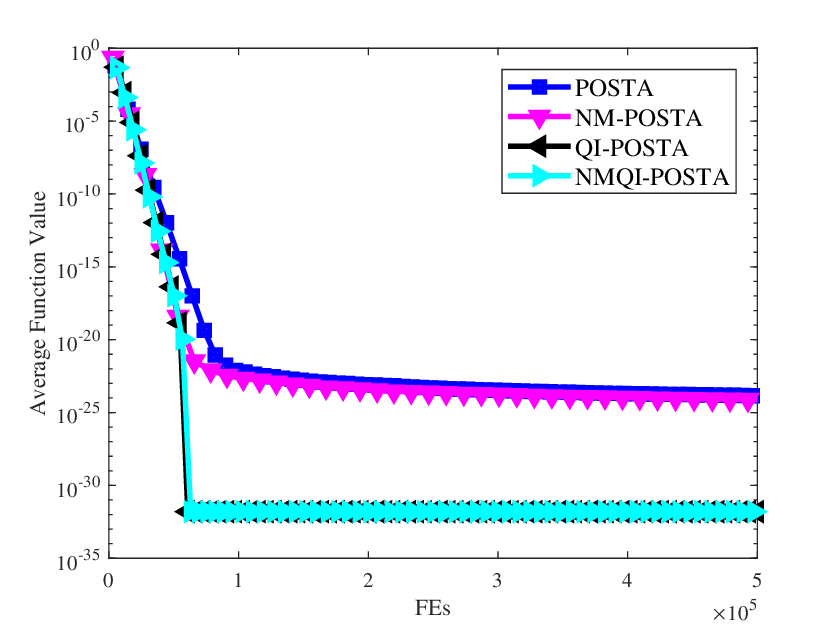}
}
\subfigure[F11] {
\label{NM2:k}
\includegraphics[width=0.2\columnwidth]{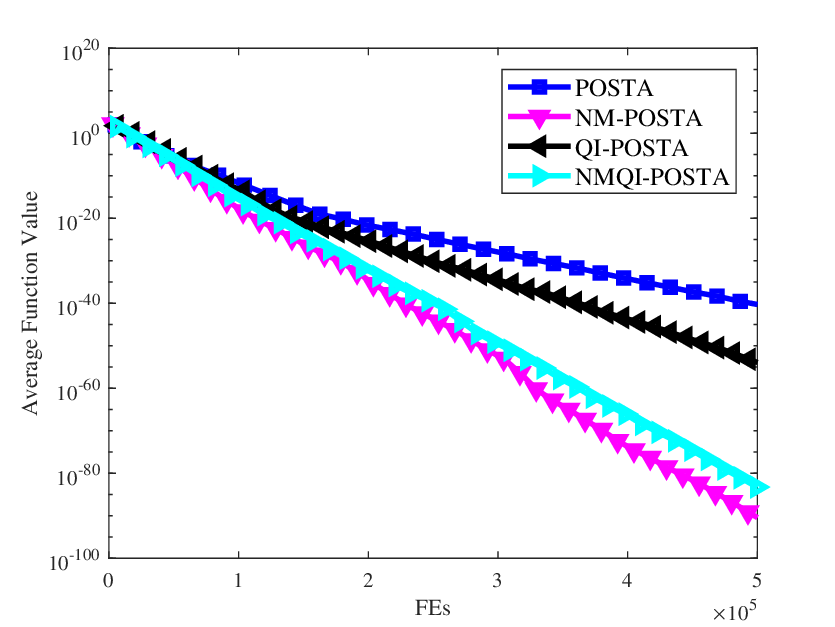}
}
\subfigure[F12] {
\label{NM2:l}
\includegraphics[width=0.2\columnwidth]{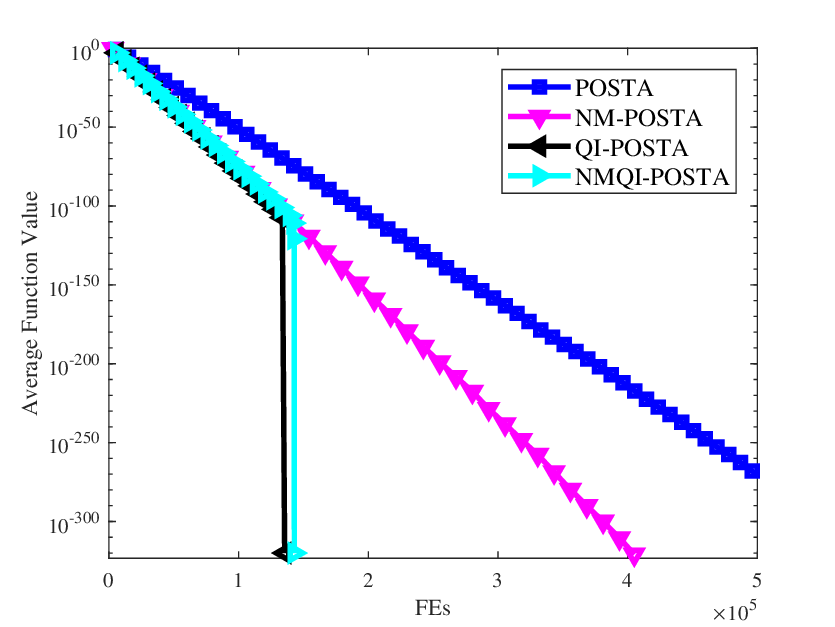}
}
\subfigure[F13] {
\label{NM2:m}
\includegraphics[width=0.2\columnwidth]{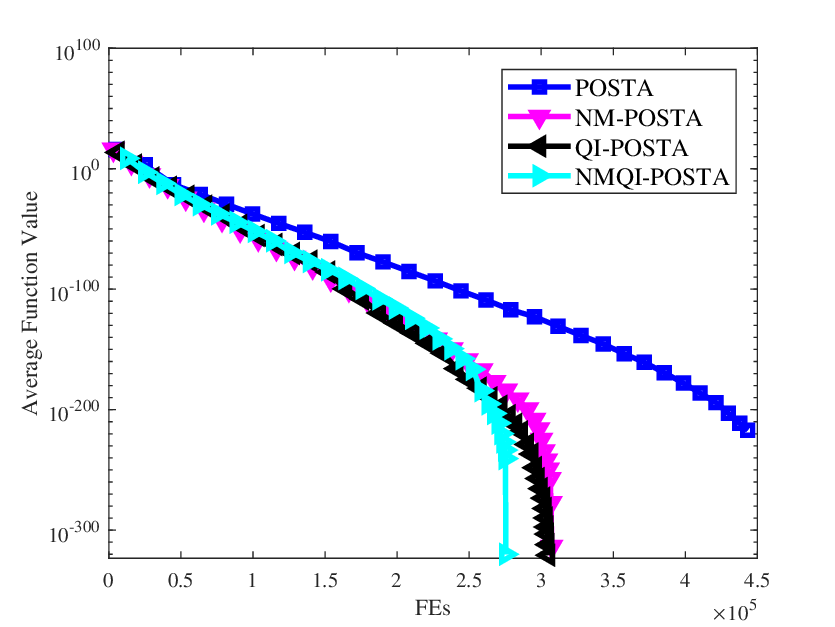}
}
\subfigure[F14] {
\label{NM2:n}
\includegraphics[width=0.2\columnwidth]{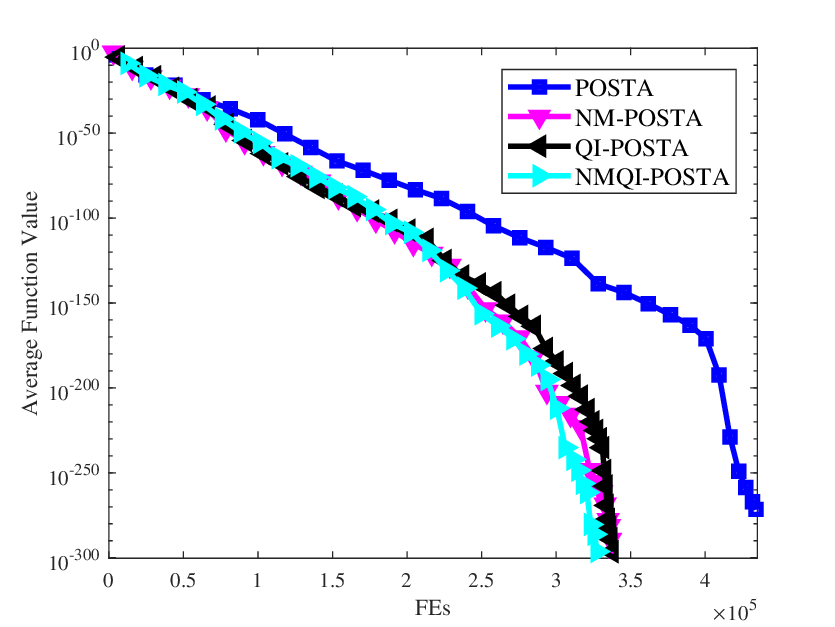}
}
\caption{Average convergence curves of the POSTA families on benchmark functions with D = 30.}
\label{NM2}
\end{figure*}

To comprehensively verify the effect of NM simplex search and QI in the proposed method, the POSTA families (POSTA, NM-POSTA, QI-POSTA and NMQI-POSTA) are tested against the benchmark functions listed in Table \ref{table3}. For each benchmark function, the number of decision variables $D$ is set to 20, 30, and 50. For each independent test, the maximum number of function evaluations is set at $5e3*D*log(D)$. In an independent test, a method is terminated if any of the following criteria is satisfied:
\begin{itemize}
\item The global minimum of the objective function is found.
\item The maximum number of function evaluations is attained.
\end{itemize}

After termination, the best solution attained and the number of function evaluations (FEs) used in this single test are recorded. The statistical results are shown in Tables \ref{table4}-\ref{table6} and the average convergence curves are given in Fig. \ref{NM2}.

The experimental results demonstrate that NM simplex search speeds up the convergence while maintaining stable, as shown in Tables \ref{table4}-\ref{table6} and Fig. \ref{NM2}. For functions like F3 and F4, POSTA suffers from a relatively low convergence speed but NM-POSTA significantly improves the speed by utilizing historical information. For function F1, F6, F9, F11, F12, F13 and F14, POSTA already gives an acceptable performance but NM-POSTA further accelerates the search process and reduces the average FEs. For multimodal functions like F2, F5, F7, F8 and F10, NM-POSTA remains stable in avoiding local minimums.

The experimental results also demonstrate that QI improves the solution accuracy and maintains stable in global exploration. As shown in Fig. \ref{NM2}, POSTA becomes very slow in the exploitation stage on F2, F5 and F10, but QI-POSTA produces solutions with much better accuracy. At the same time, for F1, F6, F9, F11, F12, F13 and F14, QI further reduces the average FEs.\ Also, the statistical results in Tables \ref{table4}-\ref{table6} reveal that the stability in multimodal functions is remained because QI is only invoked when the solutions are close to the global minimum.

Note that both NM simplex search and QI are capable to reduce the average FEs on these functions, but the approaches are different: NM simplex search utilizes historical information during the whole search, so that the whole convergence is accelerated, and therefore the average FEs is reduced. On the other hand, QI utilizes historical information in later stage search, so that an analytical solution $\mathbf{x}^{QI}$ with high accuracy is generated. $\mathbf{x}^{QI}$ usually has the accuracy that is high enough to meet the termination conditions. Therefore the later stage of the search is `skipped' and the average FEs is reduced. As shown in Fig. \ref{NM2}, the whole search is accelerated in NM-POSTA, and the later stage is `skipped' in QI-POSTA.

\begin{table*}
\centering
\setlength{\belowcaptionskip}{10pt}
\caption{ \textbf{Table 3: }Statistical results for POSTA families with D = 20.}
\renewcommand\arraystretch{1.15}
\scalebox{0.6}{
    \begin{tabular}{p{1.5cm}p{3.5cm}p{2cm}p{3.5cm}p{2cm}p{3.5cm}p{2cm}p{3.5cm}p{2cm}ccccccccc}
    \toprule
    \multirow{2}[0]{*}{function} & \multicolumn{2}{c}{POSTA} & \multicolumn{2}{c}{NM-POSTA} & \multicolumn{2}{c}{QI-POSTA} & \multicolumn{2}{c}{NMQI-POSTA} \\
    \cmidrule{2-3}\cmidrule{4-5}\cmidrule{6-7}\cmidrule{8-9}
          & Mean(Std)   & Ave FEs   & Mean(Std)   & Ave FEs   & Mean(Std)   & Ave FEs   & Mean(Std)   & Ave FEs \\
    \midrule
    F1    & 0.00E+00(0.00E+00) & 2.79E+05 & 0.00E+00(0.00E+00) & 1.93E+05 & 0.00E+00(0.00E+00) & \textbf{9.90E+04} & 0.00E+00(0.00E+00) & 1.05E+05 \\
    F2    & 1.71E-24(8.14E-25) & 3.01E+05 & 9.16E-25(5.61E-25) & 3.01E+05 & 2.36E-32(8.35E-48) & 3.01E+05 & \textbf{2.36E-32}(\textbf{8.35E-48}) & 3.01E+05 \\
    F3    & 8.22E+00(2.29E-01) & 3.03E+05 & 5.32E-19(8.97E-19) & 3.01E+05 & 5.86E+00(1.79E+00) & 3.01E+05 & \textbf{0.00E+00}(\textbf{0.00E+00}) & \textbf{2.13E+05} \\
    F4    & 3.67E-26(1.59E-25) & 3.01E+05 & 3.41E-58(8.38E-58) & 3.01E+05 & 1.65E-34(6.19E-34) & 3.01E+05 & \textbf{3.28E-59}(\textbf{1.09E-58}) & 3.01E+05 \\
    F5    & 1.19E-21(6.66E-22) & 3.01E+05 & 4.49E-22(3.85E-22) & 3.01E+05 & 0.00E+00(0.00E+00) & 5.00E+04 & \textbf{0.00E+00}(\textbf{0.00E+00}) & \textbf{4.98E+04} \\
    F6    & 0.00E+00(0.00E+00) & 2.89E+05 & 0.00E+00(0.00E+00) & 1.98E+05 & \textbf{0.00E+00}(\textbf{0.00E+00}) & \textbf{1.01E+05} & 0.00E+00(0.00E+00) & 1.05E+05 \\
    F7    & 0.00E+00(0.00E+00) & 2.77E+04 & 0.00E+00(0.00E+00) & 1.96E+04 & 0.00E+00(0.00E+00) & 2.13E+04 & \textbf{0.00E+00}(\textbf{0.00E+00})& \textbf{1.86E+04} \\
    F8    & 9.40E-03(3.01E-02) & 9.42E+04 & 1.10E-02(3.03E-02) & 7.69E+04 & 6.23E-03(1.98E-02) & 7.83E+04 & 9.91E-03(2.71E-02) & 8.64E+04 \\
    F9    & 2.50e-323(0.00E+00) & 2.90E+05 & 0.00E+00(0.00E+00) & 1.96E+05 & \textbf{0.00E+00}(\textbf{0.00E+00}) & \textbf{1.00E+05} & 0.00E+00(0.00E+00) & 1.05E+05 \\
    F10   & 1.75E-24(7.30E-25) & 3.01E+05 & 8.59E-25(4.08E-25) & 3.01E+05 & 2.36E-32(8.35E-48) & 3.01E+05 & \textbf{2.36E-32}(\textbf{8.35E-48}) & 3.01E+05 \\
    F11   & 6.69E-44(3.66E-43) & 3.01E+05 & 1.71E-85(9.31E-85) & 3.01E+05 & 1.56E-59(8.52E-59) & 3.01E+05 & \textbf{4.13E-86}(\textbf{1.63E-85}) & 3.01E+05 \\
    F12   & 4.12E-203(0.00E+00) & 3.01E+05 & 2.10E-296(0.00E+00) & 3.01E+05 & \textbf{0.00E+00}(\textbf{0.00E+00}) & \textbf{1.09E+05} & 0.00E+00(0.00E+00) & 1.13E+05 \\
    F13   & 0.00E+00(0.00E+00) & 1.74E+05 & 0.00E+00(0.00E+00) & 1.18E+05 & 0.00E+00(0.00E+00) & 1.16E+05 & 0.00E+00(0.00E+00) & \textbf{1.07E+05} \\
    F14   & 0.00E+00(0.00E+00) & 1.73E+05 & 0.00E+00(0.00E+00) & \textbf{1.30E+05} & 0.00E+00(0.00E+00) & 1.32E+05 & 0.00E+00(0.00E+00) & 1.33E+05 \\
\bottomrule
\end{tabular}}%
\label{table4}%
\end{table*}%

\begin{table*}[t]
\centering
\setlength{\belowcaptionskip}{10pt}
\caption{ \textbf{Table 4: }Statistical results for POSTA families with D = 30. }
\renewcommand\arraystretch{1.15}
\scalebox{0.6}{
\begin{tabular}{p{1.5cm}p{3.5cm}p{2cm}p{3.5cm}p{2cm}p{3.5cm}p{2cm}p{3.5cm}p{2cm}ccccccccc}
    \toprule
    \multirow{2}[0]{*}{function} & \multicolumn{2}{c}{POSTA} & \multicolumn{2}{c}{NM-POSTA} & \multicolumn{2}{c}{QI-POSTA} & \multicolumn{2}{c}{NMQI-POSTA} \\
    \cmidrule{2-3}\cmidrule{4-5}\cmidrule{6-7}\cmidrule{8-9}
          & Mean(Std)   & Ave FEs   & Mean(Std)   & Ave FEs   & Mean(Std)   & Ave FEs   & Mean(Std)   & Ave FEs \\
    \midrule
    F1    & 0.00E+00(0.00E+00) & 3.77E+05 & 0.00E+00(0.00E+00) & 2.63E+05 & 0.00E+00(0.00E+00) & \textbf{1.32E+05} & 0.00E+00(0.00E+00) & 1.40E+05 \\
    F2    & 1.49E-24(5.46E-25) & 5.12E+05 & 7.93E-25(4.16E-25) & 5.12E+05 & 1.57E-32(5.57E-48) & 5.12E+05 & \textbf{1.57E-32}(\textbf{5.57E-48}) & 5.12E+05 \\
    F3    & 2.02E+01(1.74E+01) & 5.13E+05 & 8.94E-18(2.46E-17) & 5.12E+05 & 2.01E+01(1.91E+01) & 5.12E+05 & \textbf{0.00E+00}(\textbf{0.00E+00}) & \textbf{3.18E+05} \\
    F4    & 2.01E-20(6.89E-20) & 5.12E+05 & \textbf{1.25E-74}(\textbf{3.78E-74}) & 5.12E+05 & 1.47E-25(6.19E-25) & 5.12E+05 & 6.23E-74(3.28E-73) & 5.12E+05 \\
    F5    & 1.76E-21(6.05E-22) & 5.12E+05 & 4.56E-22(3.37E-22) & 5.12E+05 & 0.00E+00(0.00E+00) & 8.04E+04 & \textbf{0.00E+00}(\textbf{0.00E+00}) & \textbf{6.94E+04} \\
    F6    & 0.00E+00(0.00E+00) & 3.90E+05 & 0.00E+00(0.00E+00) & 2.72E+05 & 0.00E+00(0.00E+00) & \textbf{1.38E+05} & 0.00E+00(0.00E+00) & 1.46E+05 \\
    F7    & 0.00E+00(0.00E+00) & 5.89E+04 & 0.00E+00(0.00E+00) & 3.25E+04 & 0.00E+00(0.00E+00) & \textbf{3.24E+04} & 0.00E+00(0.00E+00) & 3.35E+04 \\
    F8    & 0.00E+00(0.00E+00) & 4.26E+04 & 0.00E+00(0.00E+00) & 2.94E+04 & 0.00E+00(0.00E+00) & 4.85E+04 & 0.00E+00 (0.00E+00) & \textbf{2.91E+04} \\
    F9    & 0.00E+00(0.00E+00) & 3.88E+05 & 0.00E+00(0.00E+00) & 2.68E+05 & 0.00E+00(0.00E+00) & \textbf{1.39E+05} & 0.00E+00(0.00E+00) & 1.45E+05 \\
    F10   & 1.33E-24(4.26E-25) & 5.12E+05 & 6.01E-25(3.10E-25) & 5.12E+05 & 1.57E-32(5.57E-48) & 5.12E+05 & \textbf{1.57E-32}(\textbf{5.57E-48}) & 5.12E+05 \\
    F11   & 1.40E-41(4.06E-41) & 5.12E+05 & \textbf{1.20E-86}(\textbf{3.70E-86}) & 5.12E+05 & 7.10E-56(2.76E-55) & 5.12E+05 & 4.97E-86(1.77E-85) & 5.12E+05 \\
    F12   & 6.97E-280(0.00E+00) & 5.14E+05 & 0.00E+00(0.00E+00) & 4.07E+05 & 0.00E+00(0.00E+00) & \textbf{1.36E+05} & 0.00E+00(0.00E+00) & 1.43E+05 \\
    F13   & 8.17E-218(0.00E+00) & 4.43E+05 & 0.00E+00(0.00E+00) & 3.59E+05 & 0.00E+00(0.00E+00) & 3.05E+05 & 0.00E+00(0.00E+00) & \textbf{2.75E+05} \\
    F14   & 2.97E-272(0.00E+00) & 4.34E+05 & \textbf{0.00E+00}(\textbf{0.00E+00}) & \textbf{3.25E+05} & 1.48E-299(0.00E+00) & 3.38E+05 & 0.00E+00(0.00E+00) & 3.30E+05 \\
\bottomrule
\end{tabular}}%
\label{table5}%
\end{table*}%

\begin{table*}[h]
\centering
\setlength{\belowcaptionskip}{10pt}
\caption{ \textbf{Table 5: }Statistical results for POSTA families with D = 50. }
\renewcommand\arraystretch{1.15}
\scalebox{0.6}{
    \begin{tabular}{p{1.5cm}p{3.5cm}p{2cm}p{3.5cm}p{2cm}p{3.5cm}p{2cm}p{3.5cm}p{2cm}ccccccccc}
    \toprule
    \multirow{2}[0]{*}{function} & \multicolumn{2}{c}{POSTA} & \multicolumn{2}{c}{NM-POSTA} & \multicolumn{2}{c}{QI-POSTA} & \multicolumn{2}{c}{NMQI-POSTA} \\
  \cmidrule{2-3}\cmidrule{4-5}\cmidrule{6-7}\cmidrule{8-9}
          & Mean(Std)   & Ave FEs   & Mean(Std)   & Ave FEs   & Mean(Std)   & Ave FEs   & Mean(Std)   & Ave FEs \\
    \midrule
    F1    & 0.00E+00(0.00E+00) & 5.84E+05 & 0.00E+00(0.00E+00) & 4.23E+05 & 0.00E+00(0.00E+00) & \textbf{2.17E+05} & 0.00E+00(0.00E+00) & 2.24E+05 \\
    F2    & 9.85E-25(2.38E-25) & 9.80E+05 & 4.59E-25(1.85E-25) & 9.80E+05 & 9.42E-33(2.78E-48) & 9.80E+05 & \textbf{9.42E-33}(\textbf{2.78E-48}) & 9.80E+05 \\
    F3    & 4.44E+01(2.31E+01) & 9.80E+05 & 1.42E-16(1.74E-16) & 9.80E+05 & 3.40E+01(1.94E+01) & 9.80E+05 & \textbf{0.00E+00}(\textbf{0.00E+00}) & \textbf{6.61E+05} \\
    F4    & 1.18E-12(2.07E-12) & 9.80E+05 & \textbf{1.38E-87}(\textbf{7.50E-87}) & 9.80E+05 & 7.79E-17(1.43E-16) & 9.80E+05 & 7.21E-86(3.94E-85) & 9.80E+05 \\
    F5    & 2.89E-21(9.30E-22) & 9.80E+05 & 2.10E-21(1.29E-21) & 9.80E+05 & 0.00E+00(0.00E+00) & 1.41E+05 & 0.00E+00(0.00E+00) & \textbf{9.49E+04} \\
    F6    & 0.00E+00(0.00E+00) & 6.24E+05 & 0.00E+00(0.00E+00) & 4.61E+05 & 0.00E+00(0.00E+00) & \textbf{2.42E+05} & 0.00E+00(0.00E+00) & 2.63E+05 \\
    F7    & 3.60E-14(7.08E-14) & 4.08E+05 & 4.17E-14(8.70E-14) & 3.52E+05 & 3.79E-14(7.51E-14) & 3.25E+05 & 5.12E-14(8.22E-14) & 4.61E+05 \\
    F8    & 0.00E+00(0.00E+00) & 7.33E+04 & 0.00E+00(0.00E+00) & \textbf{5.51E+04} & 0.00E+00(0.00E+00) & 6.54E+04 & 0.00E+00(0.00E+00) & 5.93E+04 \\
    F9    & 0.00E+00(0.00E+00) & 6.38E+05 & 0.00E+00(0.00E+00) & 4.43E+05 & 0.00E+00(0.00E+00) & \textbf{2.32E+05} & 0.00E+00(0.00E+00) & 2.47E+05 \\
    F10   & 9.15E-25(2.12E-25) & 9.80E+05 & 5.12E-25(2.03E-25) & 9.80E+05 & 9.42E-33(2.78E-48) & 9.80E+05 & \textbf{9.42E-33}(\textbf{2.78E-48}) & 9.80E+05 \\
    F11   & 1.54E-38(3.96E-38) & 9.80E+05 & 2.53E-105(8.38E-105) & 9.80E+05 & 6.76E-50(3.35E-49) & 9.80E+05 & \textbf{1.86E-106}(\textbf{4.36E-106}) & 9.80E+05 \\
    F12   & 0.00E+00(0.00E+00) & 7.47E+05 & 0.00E+00(0.00E+00) & 5.30E+05 & 0.00E+00(0.00E+00) & \textbf{1.74E+05} & 0.00E+00(0.00E+00) & 1.85E+05 \\
    F13   & 7.33E-276(0.00E+00) & 9.79E+05 & 0.00E+00(0.00E+00) & 7.91E+05 & 0.00E+00(0.00E+00) & 6.07E+05 & \textbf{0.00E+00}(\textbf{0.00E+00}) & \textbf{6.03E+05} \\
    F14   & 6.39E-191(0.00E+00) & 9.80E+05 & \textbf{2.10E-267}(\textbf{0.00E+00}) & 9.80E+05 & 7.44E-259(0.00E+00) & 9.80E+05 & 2.26E-266(0.00E+00) & 9.80E+05 \\
\bottomrule
\end{tabular}}%
\label{table6}%
\end{table*}%

As mentioned above, NM simplex search is applied during the whole search. But QI is only invoked in later stage. As shown in Fig. \ref{NM2:c}, QI-POSTA shows no advantage over POSTA on F3 because POSTA has a low convergence speed on F3 and is not able to reach the neighborhood of the global minimum. But in NMQI-POSTA, QI is quite effective because NM-POSTA is able to reach the neighborhood of the global minimum.

From the above discussion, it can be observed that when POSTA is enhanced with both NM simplex search and QI, the merits of the three distinct methods are combined. The convergence curves and the statistical results show that NMQI-POSTA enjoys the merits of the three: global exploration capacity of POSTA, fast convergence speed of NM simplex search and deep exploitation capacity of QI.

~\\
\textbf{4.3. Comparison with other meteheuristic methods}
~\\

To examine the effectiveness of the enhanced POSTA, several population-based metaheuristic methods are used in this section, including ABC \cite{ABC2007, ABC2020spam}, SaDE \cite{SaDE2008,sade2021self}, GWO \cite{GWO2014,GWO2021swarm} and CLPSO \cite{CLPSO2006,CLPSO2020swarm}. For each benchmark function, the number of decision variables $D$ is set to 20, 30, and 50. For each independent test, the maximum function evaluations is set at $5e3*D*log(D)$.

\begin{table}
\centering
\setlength{\belowcaptionskip}{10pt}
\caption{ \textbf{Table 6: }Statistical results for different metaheuristic methods with D = 20. }
\renewcommand\arraystretch{1.15}
\scalebox{0.6}{
    \begin{tabular}{ccccccccccccccc}
    \toprule
    \multirow{2}[0]{*}{function} & \multicolumn{1}{c}{ABC} & \multicolumn{1}{c}{GWO} & \multicolumn{1}{c}{SaDE} & \multicolumn{1}{c}{CLPSO} & \multicolumn{1}{c}{NMQI-POSTA} \\
    \cmidrule{2-6}
          & Mean(Std)  & Mean(Std)  & Mean(Std)   & Mean(Std)  & Mean(Std) \\
    \midrule
    F1    & 3.42E+03(2.99E+03)- & 0.00E+00(0.00E+00)$\approx$ & 3.93E-201(0.00E+00)-     & 1.20E-33(1.73E-33)-     & \textbf{0.00E+00}(\textbf{0.00E+00}) \\
    F2    & 2.64E-02(3.75E-02)- & 1.70E-02(7.90E-03)-     & 2.36E-32(8.35E-48)$\approx$     & 2.36E-32(8.35E-48$\approx$     & \textbf{2.36E-32}(\textbf{8.35E-48}) \\
    F3    & 4.65E+02(2.64E+02)- & 1.63E+01(6.38E-01)-     & 9.81E+00(2.84E+00)-     & 2.19E+00(2.30E+00)-     & \textbf{0.00E+00}(\textbf{0.00E+00}) \\
    F4    & 1.02E+04(2.06E+03)- & \textbf{1.53E-301}(\textbf{0.00E+00})+ & 3.22E-20(7.57E-20)-     & 7.86E+00(2.53E+00)-     & 3.28E-59(1.09E-58) \\
    F5    & 1.12E+01(1.01E+01)- & 3.84E+00(1.21E+00)-     & 1.58E-30(7.61E-31)-     & 3.25E-05(7.76E-06)-     & \textbf{0.00E+00}(\textbf{0.00E+00}) \\
    F6    & 4.26E-01(3.24E-01)- & 0.00E+00(0.00E+00)$\approx$ & 1.73E-197(0.00E+00)-     & 2.90E-37(2.62E-37)-     & 0.00E+00(0.00E+00) \\
    F7    & 1.09E+01(2.10E+00)- & 0.00E+00(0.00E+00)$\approx$ & 0.00E+00(0.00E+00)$\approx$     & 0.00E+00(0.00E+00)$\approx$     & \textbf{0.00E+00}(\textbf{0.00E+00}) \\
    F8    & 5.92E-01(1.51E-01)- & \textbf{0.00E+00}(\textbf{0.00E+00})+     & 5.75E-04(2.21E-03)$\approx$     & 2.59E-17(9.07E-17)$\approx$     & 9.91E-03(2.71E-02) \\
    F9    & 4.45E-02(2.51E-02)- & 0.00E+00(0.00E+00)$\approx$ & 2.38E-199(0.00E+00)-     & 2.16E-38(1.64E-38)-     & \textbf{0.00E+00}(\textbf{0.00E+00}) \\
    F10   & 1.00E-04(1.22E-04)- & 1.82E-02(1.03E-02)- & 2.36E-32(8.35E-48)$\approx$     & 2.36E-32(8.35E-48)$\approx$ & \textbf{2.36E-32}(\textbf{8.35E-48}) \\
    F11   & 1.70E+02(2.49E+01)- & \textbf{0.00E+00}(\textbf{0.00E+00})+ & 1.26E-33(3.95E-33)-     & 5.69E-04(1.78E-04)- & 4.13E-86(1.63E-85) \\
    F12   & 1.54E-01(3.88E-02)- & 0.00E+00(0.00E+00)$\approx$     & 2.18E-107(2.80E-107)-     & 1.99E-22(1.04E-22)-     & \textbf{0.00E+00}(\textbf{0.00E+00}) \\
    F13   & 3.02E+04(9.22E+04)- & 0.00E+00(0.00E+00)$\approx$     & 7.55e-311(0.00E+00)-     & 2.08E-66(3.88E-66)-     & \textbf{0.00E+00}(\textbf{0.00E+00}) \\
    F14   & 7.69E-14(2.55E-13)- & 0.00E+00(0.00E+00)$\approx$     & 0.00E+00(0.00E+00)$\approx$     & 1.34E-78(2.15E-78)-     & \textbf{0.00E+00}(\textbf{0.00E+00}) \\
    +/$\approx$/-& 0/0/14 & 3/7/4 & 0/5/9 & 0/4/10 &-/-/- \\
\bottomrule
\end{tabular}}%
\label{table7}%
\end{table}%

\begin{figure*}
\centering
\subfigure[F1] {
\label{ALLMETHOD1:a}
\includegraphics[width=0.2\columnwidth]{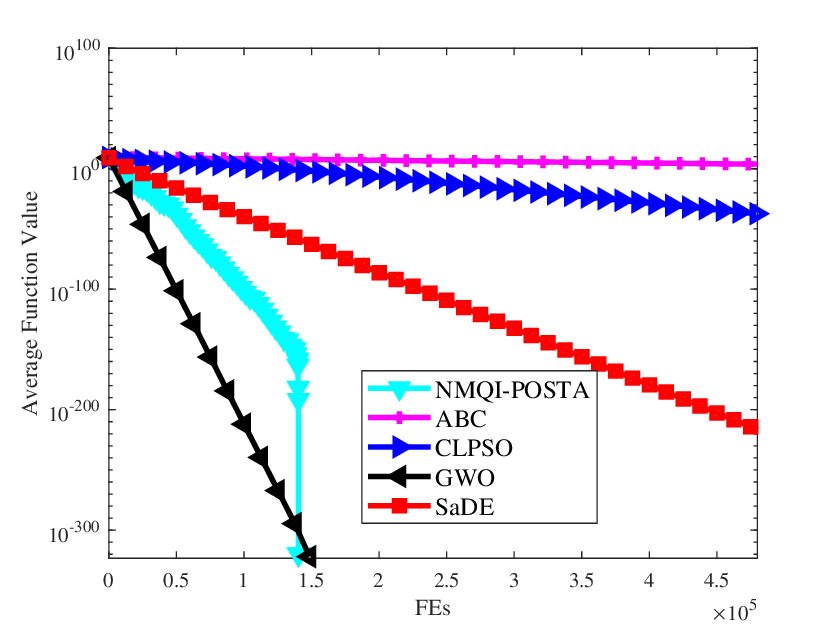}
}
\subfigure[F2] {
\label{ALLMETHOD1:b}
\includegraphics[width=0.2\columnwidth]{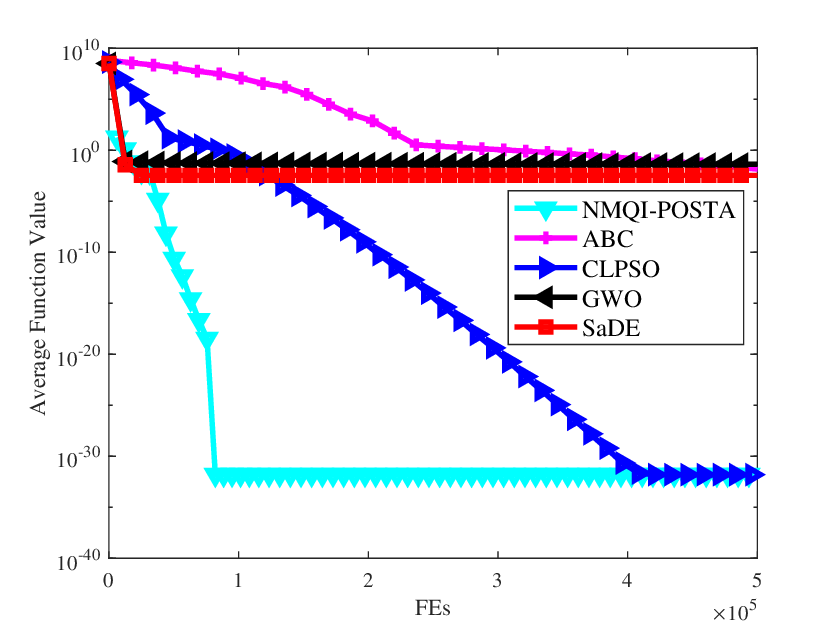}
}
\subfigure[F3] {
\label{ALLMETHOD1:c}
\includegraphics[width=0.2\columnwidth]{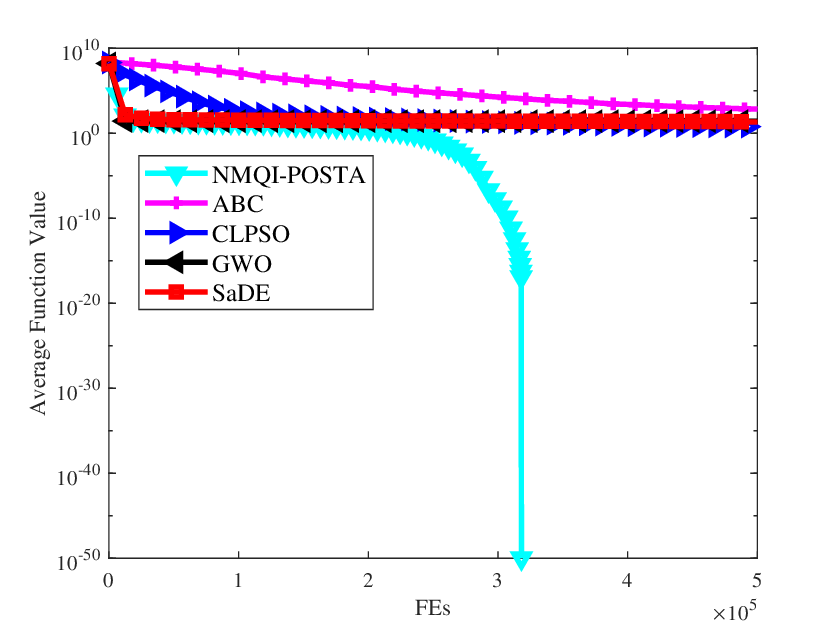}
}
\subfigure[F4] {
 \label{ALLMETHOD1:d}
\includegraphics[width=0.2\columnwidth]{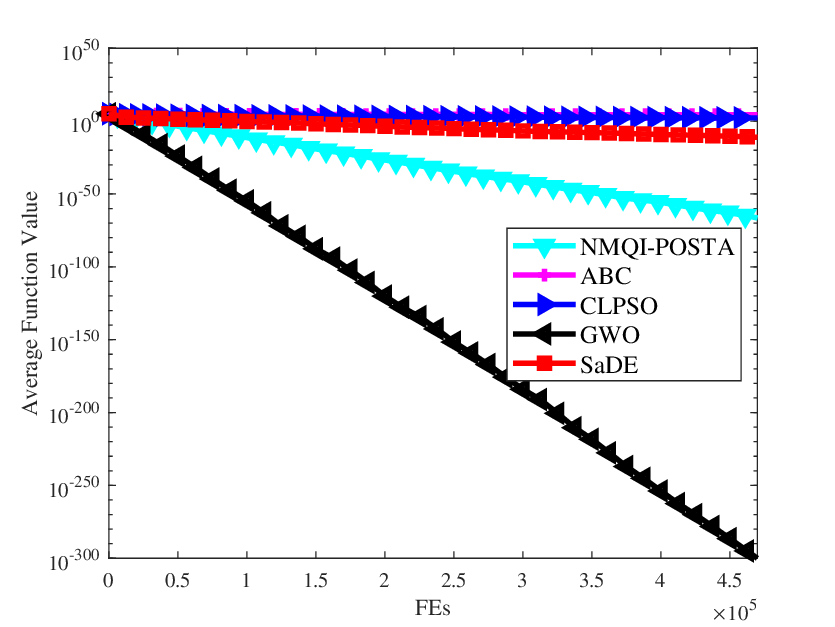}
}
\subfigure[F5] {
\label{ALLMETHOD1:e}
\includegraphics[width=0.2\columnwidth]{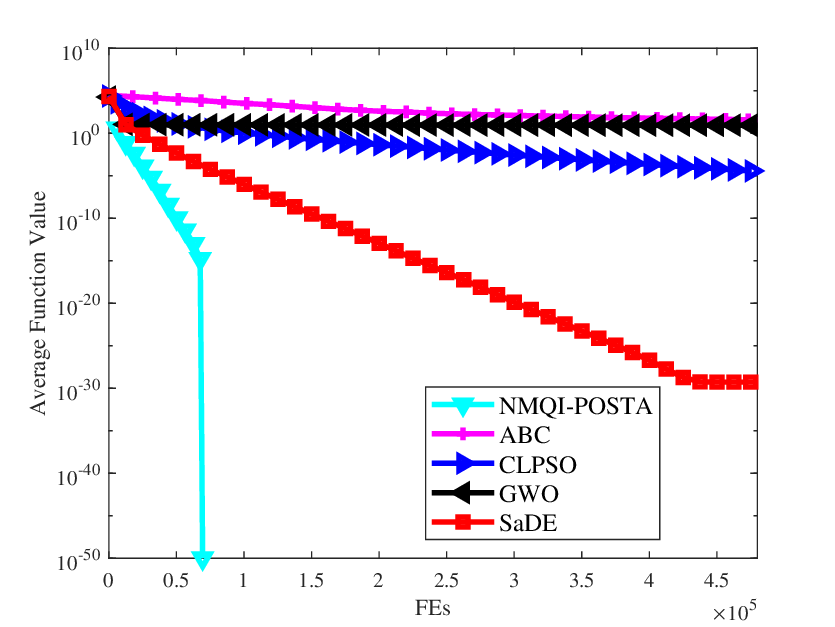}
}
\subfigure[F6] {
\label{ALLMETHOD1:f}
\includegraphics[width=0.2\columnwidth]{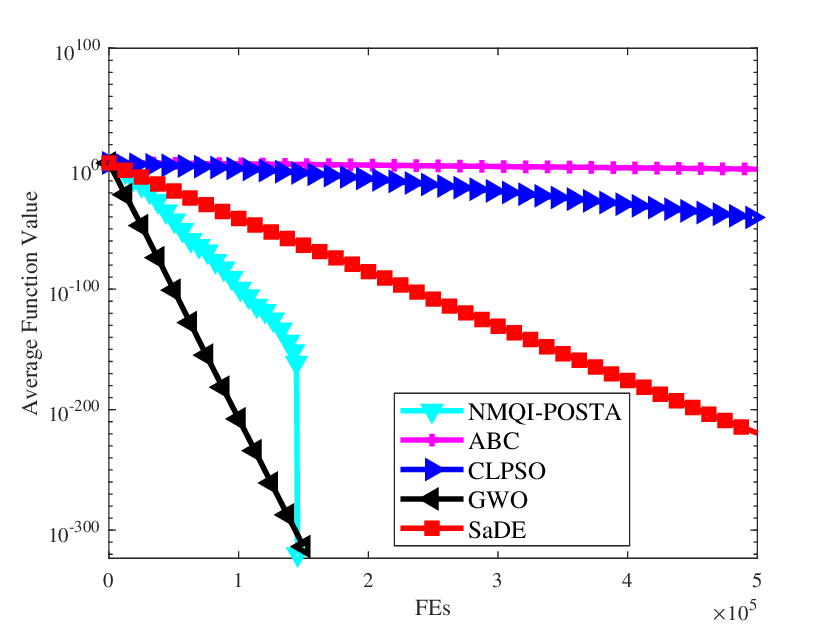}
}
\subfigure[F7] {
\label{ALLMETHOD1:g}
\includegraphics[width=0.2\columnwidth]{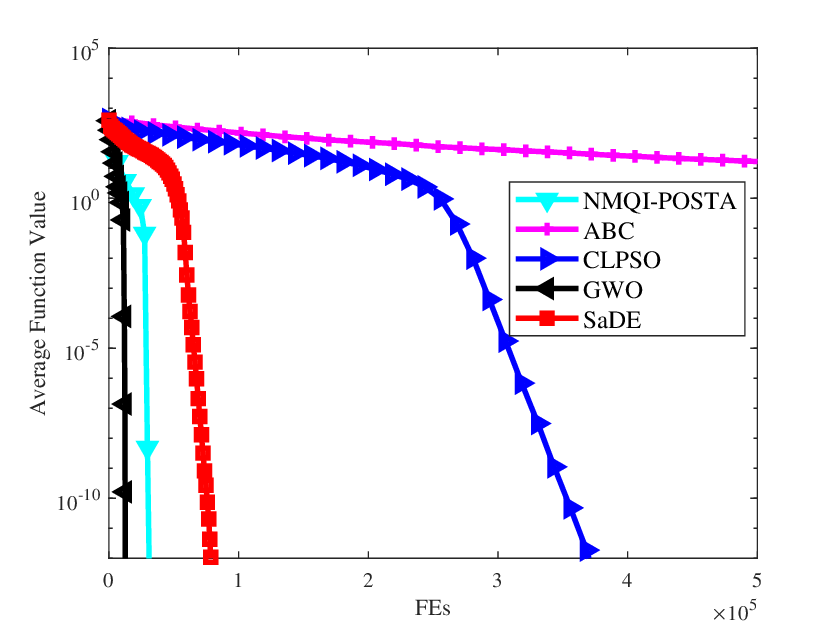}
}
\subfigure[F8] {
\label{ALLMETHOD1:h}
\includegraphics[width=0.2\columnwidth]{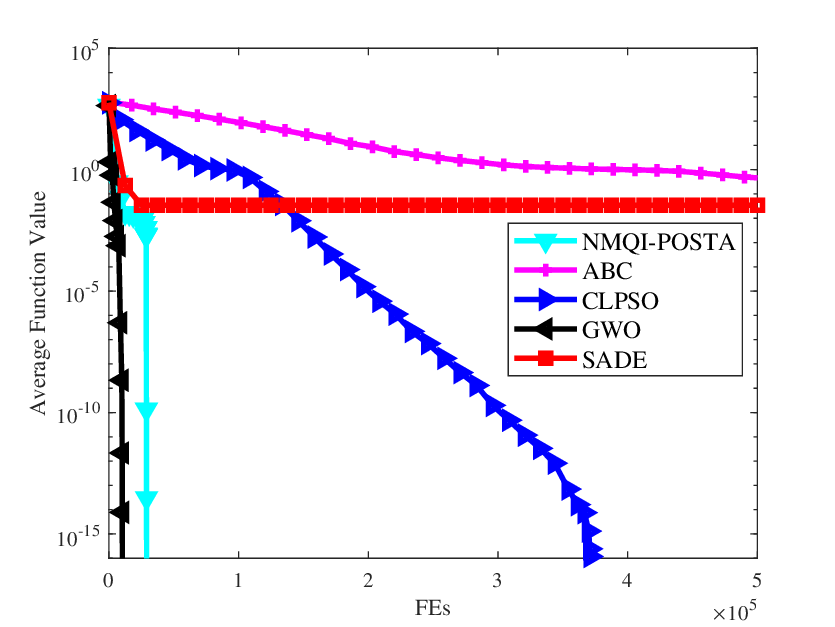}
}
\subfigure[F9] {
\label{ALLMETHOD1:i}
\includegraphics[width=0.2\columnwidth]{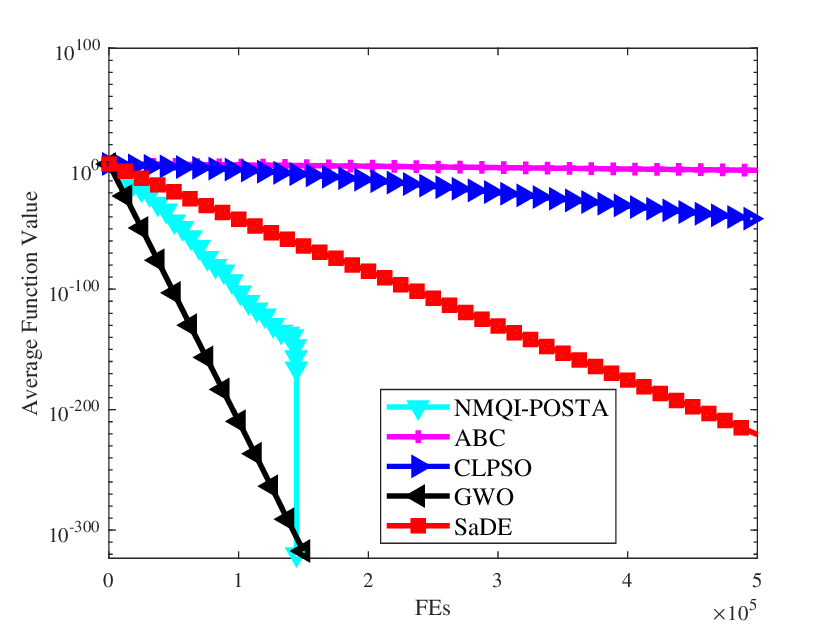}
}
\subfigure[F10] {
\label{ALLMETHOD1:j}
\includegraphics[width=0.2\columnwidth]{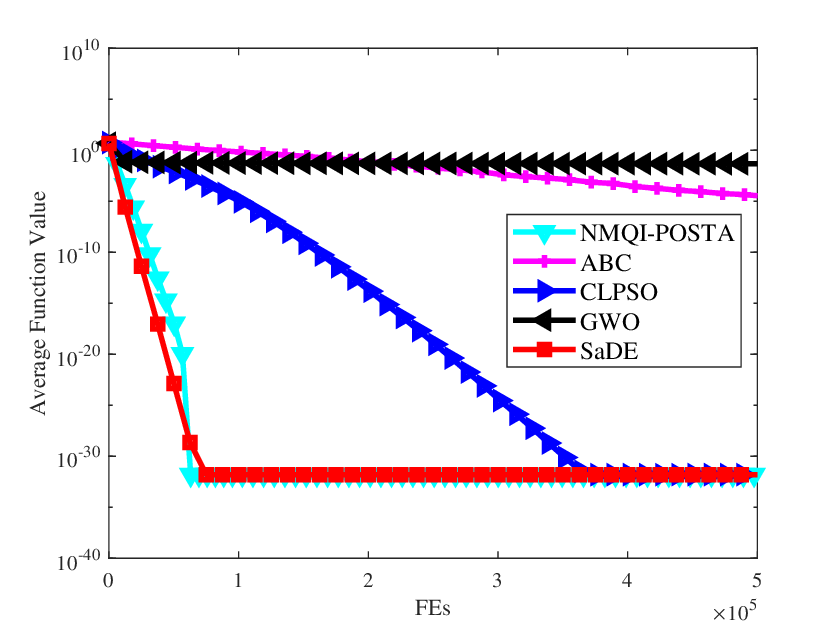}
}
\subfigure[F11] {
\label{ALLMETHOD1:k}
\includegraphics[width=0.2\columnwidth]{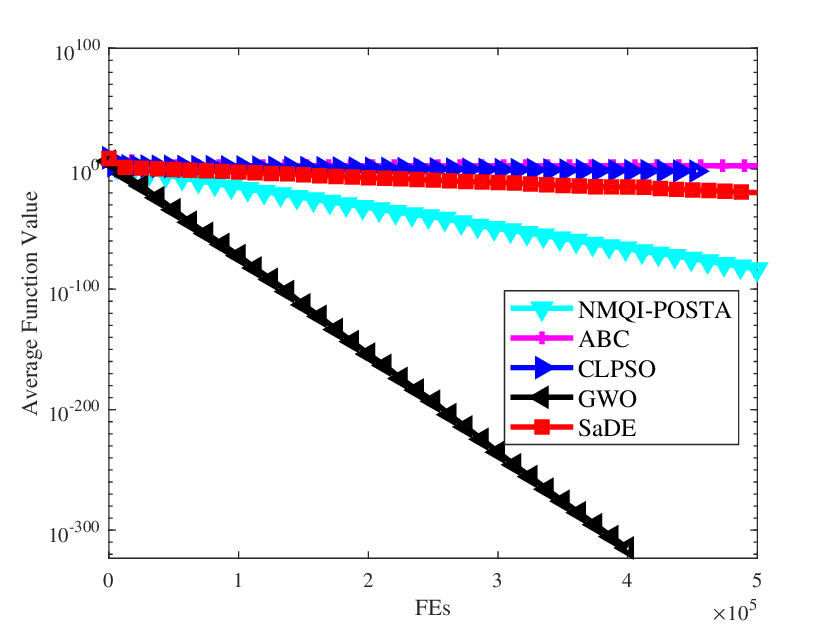}
}
\subfigure[F12] {
\label{ALLMETHOD1:l}
\includegraphics[width=0.2\columnwidth]{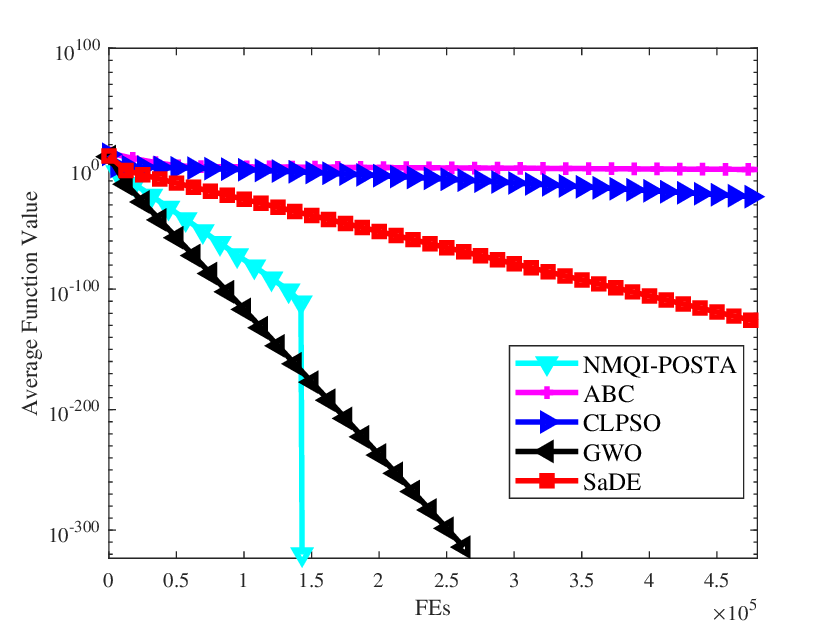}
}
\subfigure[F13] {
\label{ALLMETHOD1:m}
\includegraphics[width=0.2\columnwidth]{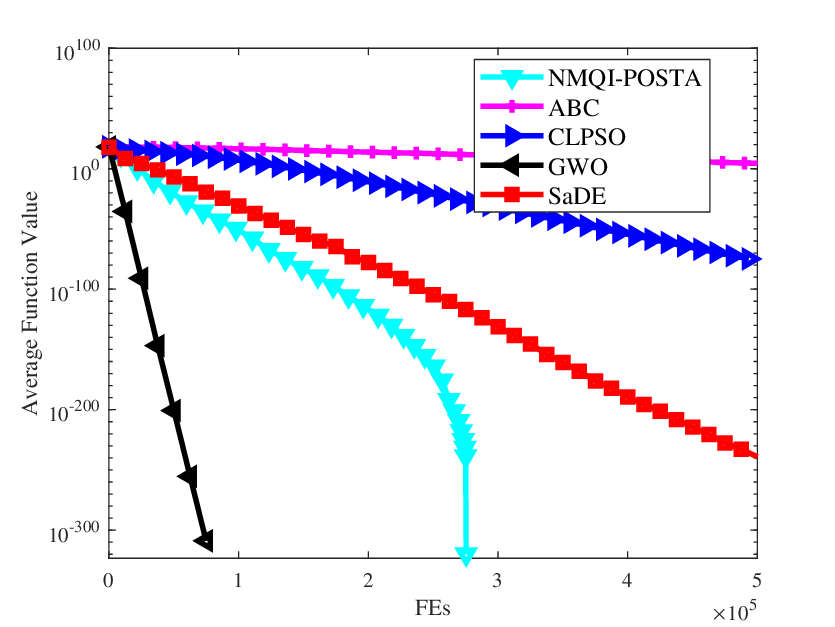}
}
\subfigure[F14] {
\label{ALLMETHOD1:n}
\includegraphics[width=0.2\columnwidth]{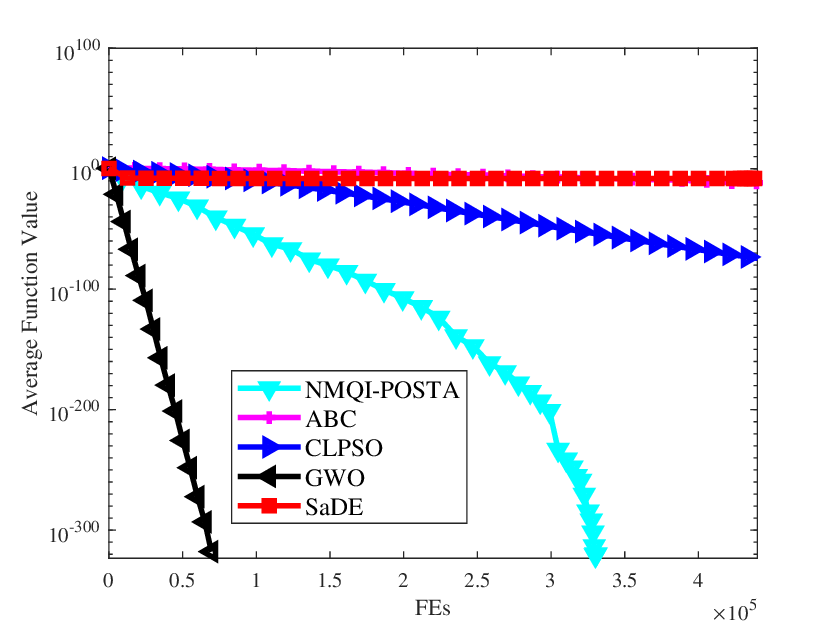}
}
\caption{Average convergence curves of different metaheuristic methods on benchmark functions with D = 30. }
\label{ALLMETHOD1}
\end{figure*}

The computational results are listed in Tables \ref{table7}-\ref{table9} and the average convergence curves of all test functions with $D = 30$ are given in Fig. \ref{ALLMETHOD1}. In Tables \ref{table7}-\ref{table9}, the symbols -, +, $\approx$ indicate the results of Wilcoxon rank sum test \cite{rst2021swarm1,rst2021swarm2} with a significance level of 0.05. Correspondingly, the symbol - indicates that the performance is significantly worse than the proposed method, the symbol + indicates that the performance is significantly better than the proposed method, and the symbol $\approx$ indicates that the performance difference is not statistically significant compared to the proposed method.

As shown in the experiment results, the enhanced POSTA outperforms CLPSO, SaDE and ABC in most of the cases. Despite GWO gives better performance on functions like F4 and F11, it failed to remain stable in finding acceptable solutions on F2, F3, F5 and F10. On the contrary, the enhanced POSTA shows incomparably convergence speed and solution accuracy on F3 and F5, and shows more robustness on multimodal functions like F2 and F10. Therefore, the proposed method is more effective and robust than these competitive metaheuristic methods.

\begin{table}
\centering
\setlength{\belowcaptionskip}{10pt}
\caption{ \textbf{Table 7: }Statistical results for different metaheuristic methods with D = 30. }
\renewcommand\arraystretch{1.15}
\scalebox{0.6}{
    \begin{tabular}{ccccccccccccccc}
    \toprule
    \multirow{2}[0]{*}{function} & \multicolumn{1}{c}{ABC} & \multicolumn{1}{c}{GWO} & \multicolumn{1}{c}{SaDE} & \multicolumn{1}{c}{CLPSO} & \multicolumn{1}{c}{NMQI-POSTA} \\
    \cmidrule{2-6}
          & Mean(Std)  & Mean(Std)  & Mean(Std)   & Mean(Std)  & Mean(Std) \\
    \midrule
    F1    & 2.53E+03(1.88E+03)-     & 0.00E+00(0.00E+00)$\approx$     & 5.72E-231(0.00E+00)-     & 3.95E-39(3.51E-39)-     & \textbf{0.00E+00}(\textbf{0.00E+00}) \\
    F2    & 1.07E-02(1.13E-02)-     & 4.27E-02(1.38E-02)-     & 1.38E-02(7.57E-02)-     & 1.57E-32(5.57E-48)$\approx$     & \textbf{1.57E-32}(\textbf{5.57E-48}) \\
    F3    & 6.31E+02(3.37E+02)-     & 2.63E+01(6.65E-01)-     & 1.87E+01(1.55E+01)-     & 5.87E+00(1.39E+01)-     & \textbf{0.00E+00}(\textbf{0.00E+00}) \\
    F4    & 2.31E+04(2.93E+03)-     & \textbf{0.00E+00}(\textbf{0.00E+00})+     & 2.35E-13(1.03E-12)-     & 1.52E+02(4.31E+01)-     & 6.23E-74(3.28E-73) \\
    F5    & 3.05E+01(1.45E+01)-     & 8.80E+00(8.62E-01)-     & 5.27E-30(1.83E-30)-     & 2.95E-05(6.95E-06)-     & \textbf{0.00E+00}(\textbf{0.00E+00}) \\
    F6    & 3.00E-01(2.01E-01)-     & 0.00E+00(0.00E+00)$\approx$     & 1.19E-224(0.00E+00)-     & 2.61E-42(1.72E-42)-     & \textbf{0.00E+00}(\textbf{0.00E+00}) \\
    F7    & 1.60E+01(2.54E+00)-     & 0.00E+00(0.00E+00)$\approx$     & 0.00E+00(0.00E+00)$\approx$     & 0.00E+00(0.00E+00)$\approx$     & \textbf{0.00E+00}(\textbf{0.00E+00}) \\
    F8    & 4.27E-01(1.21E-01)-     & 0.00E+00(0.00E+00)$\approx$     & 3.86E-03(8.49E-03)-     & 0.00E+00(0.00E+00)$\approx$     & \textbf{0.00E+00}(\textbf{0.00E+00}) \\
    F9    & 4.02E-02(2.32E-02)-     & 0.00E+00(0.00E+00)$\approx$     & 5.24E-226(0.00E+00)-     & 2.98E-43(2.53E-43)-     & \textbf{0.00E+00}(\textbf{0.00E+00}) \\
    F10   & 2.55E-05(2.79E-05)-     & 4.66E-02(1.63E-02)-     & 1.57E-32(5.57E-48)$\approx$     & 1.57E-32(5.57E-48)$\approx$     & \textbf{1.57E-32}(\textbf{5.57E-48}) \\
    F11   & 3.08E+02(3.11E+01)-     & \textbf{0.00E+00}(\textbf{0.00E+00})+     & 5.90E-21(2.09E-20)-     & 8.38E-03(3.13E-03)-     & 4.97E-86(1.77E-85) \\
    F12   & 1.46E-01(3.54E-02)-     & 0.00E+00(0.00E+00)$\approx$     & 9.29E-136(1.74E-135)-     & 2.01E-25(8.32E-26)-     & \textbf{0.00E+00}(\textbf{0.00E+00}) \\
    F13   & 1.28E+04(4.48E+04)-     & 0.00E+00(0.00E+00)$\approx$     & 2.38E-245(0.00E+00)-     & 1.71E-77(6.14E-77)-     & \textbf{0.00E+00}(\textbf{0.00E+00}) \\
    F14   & 1.07E-14(2.08E-14)-     & 0.00E+00(0.00E+00)$\approx$     & 8.05E-09(2.89E-08)-     & 2.88E-86(3.84E-86)-     & \textbf{0.00E+00}(\textbf{0.00E+00}) \\
    +/$\approx$/-& 0/0/14 & 2/8/4 & 0/2/12 & 0/4/10 &-/-/- \\
\bottomrule
\end{tabular}}%
\label{table8}%
\end{table}%

\begin{table}[h]
\centering
\setlength{\belowcaptionskip}{10pt}
\caption{ \textbf{Table 8: }Statistical results for different metaheuristic methods with D = 50. }
\renewcommand\arraystretch{1.15}
\scalebox{0.6}{
    \begin{tabular}{ccccccccccccccc}
    \toprule
    \multirow{2}[0]{*}{function} & \multicolumn{1}{c}{ABC} & \multicolumn{1}{c}{GWO} & \multicolumn{1}{c}{SaDE} & \multicolumn{1}{c}{CLPSO} & \multicolumn{1}{c}{NMQI-POSTA} \\
    \cmidrule{2-6}
          & Mean(Std)  & Mean(Std)  & Mean(Std)   & Mean(Std)  & Mean(Std) \\
    \midrule
    F1    & 1.31E+03(8.41E+02)-     & 0.00E+00(0.00E+00)$\approx$     & 3.37E-237(0.00E+00)-     & 1.33E-46(6.70E-47)-     & \textbf{0.00E+00}(\textbf{0.00E+00}) \\
    F2    & 2.82E-03(1.86E-03)-     & 1.04E-01(2.34E-02)-     & 1.45E-02(3.13E-02)-     & 9.42E-33(2.78E-48)$\approx$     & \textbf{9.42E-33}(\textbf{2.78E-48}) \\
    F3    & 6.13E+02(1.92E+02)-     & 4.63E+01(7.82E-01)-     & 4.42E+01(2.78E+01)-     & 1.32E+01(1.84E+01)-     & \textbf{0.00E+00}(\textbf{0.00E+00}) \\
    F4    & 5.49E+04(5.14E+03)-     & \textbf{0.00E+00}(\textbf{0.00E+00})+     & 5.24E-08(8.90E-08)-     & 2.55E+03(4.70E+02)-     & 7.21E-86(3.94E-85) \\
    F5    & 6.13E+01(2.15E+01)-     & 2.07E+01(1.71E+00)-     & 2.92E-29(1.32E-29)-     & 1.56E-05(2.08E-06)-     & \textbf{0.00E+00}(\textbf{0.00E+00}) \\
    F6    & 1.25E-01(6.41E-02)-     & 0.00E+00(0.00E+00)$\approx$     & 6.84E-235(0.00E+00)-     & 1.10E-49(5.63E-50)-     & \textbf{0.00E+00}(\textbf{0.00E+00}) \\
    F7    & 2.29E+01(2.47E+00)-     & \textbf{0.00E+00}(\textbf{0.00E+00})+     & 3.65E-01(6.12E-01)$\approx$     & 0.00E+0(0.00E+00)+     & 5.12E-14(8.22E-14) \\
    F8    & 2.27E-01(1.01E-01-     & 0.00E+00(0.00E+00)$\approx$     & 9.30E-03(1.95E-02)-     & 0.00E+00(0.00E+00)$\approx$     & \textbf{0.00E+00}(\textbf{0.00E+00}) \\
    F9    & 2.68E-02(1.51E-02)-     & 0.00E+00(0.00E+00)$\approx$     & 2.27E-233(0.00E+00)-     & 2.52E-50(1.56E-50)-     & \textbf{0.00E+00}(\textbf{0.00E+00}) \\
    F10   & 8.33E-06(9.02E-06)-     & 9.69E-02(2.13E-02)-     & 9.42E-33(2.78E-48)$\approx$    & 9.42E-33(2.78E-48)$\approx$     & \textbf{9.42E-33}(\textbf{2.78E-48}) \\
    F11   & 6.29E+02(5.20E+01)-     & \textbf{0.00E+00}(\textbf{0.00E+00})+     & 3.14E-09(1.69E-08)-     & 1.99E-01(6.82E-02)-     & 1.86E-106(4.36E-106) \\
    F12   & 1.11E-01(2.53E-02)-     & 0.00E+00(0.00E+00)$\approx$     & 4.89E-158(1.34E-157)-     & 8.50E-30(3.23E-30)-     & \textbf{0.00E+00}(\textbf{0.00E+00}) \\
    F13   & 6.76E+02(1.65E+03)-     & 0.00E+00(0.00E+00)$\approx$    & 1.55E-200(0.00E+00)-     & 7.03E-94(8.26E-94)-     & \textbf{0.00E+00}(\textbf{0.00E+00}) \\
    F14   & 2.08E-15(2.98E-15)-     & \textbf{0.00E+00}(\textbf{0.00E+00})+     & 6.66E-08(1.22E-07)$\approx$     & 1.34E-96(1.84E-96)-     & 2.26E-266(0.00E+00) \\
    +/$\approx$/-& 0/0/14 & 4/6/4 & 0/3/11 & 0/3/11 &-/-/- \\
\bottomrule
\end{tabular}}%
\label{table9}%
\end{table}

\section{Conclusion}
In this paper, an enhanced POSTA based on Nelder-Mead simplex search and QI is proposed. In the enhanced POSTA, both NM simplex search and QI are applied to utilize the historical information. Specifically, NM simplex search is applied in the whole search process to speed up the convergence, and QI is applied in the later stage of the search to improve the solution accuracy.
The proposed method enjoys the merits of the three methods: global exploration capacity of POSTA, fast convergence speed of NM simplex search and strong exploitation ability of QI. To demonstrate the effectiveness of the proposed method, the enhanced POSTA is tested against 14 benchmark functions in 20-D, 30-D and 50-D space. An experimental comparison with other well-known metaheuristic methods also demonstrates the robustness and the effectiveness of the enhanced POSTA.

In the proposed method, a collection strategy is used to control the selection of vertices in the initial simplex set. However, the selection of vertices requires further study. In our future study, the quantitative properties of the vertices will be considered in the collection of historical information. It would also be interesting to apply the proposed historical information mechanism to other metaheuristic methods. Besides, other approaches of utilizing the historical information will be considered as well.

\bibliographystyle{unsrt}
\bibliography{cas-refs}
\end{document}